\documentclass[10pt,reqno,oneside]{amsart}
 \usepackage{verbatim,amsmath,amssymb,cite,xspace}
\usepackage{color,cite,graphicx}

\theoremstyle{plain}
\newtheorem{theorem}{Theorem}[section]
\newtheorem{claim}{Claim}[section]
\newtheorem{lemma}{Lemma}[section]
\newtheorem{corollary}{Corollary}[section]
\theoremstyle{definition}
\newtheorem{definition}{Definition}[section]

\theoremstyle{remark}

\usepackage{color}

\usepackage{amssymb}
\usepackage{amssymb, latexsym}
\usepackage{amsmath}
\usepackage{amscd}
\usepackage{enumerate}
\usepackage{amsfonts}
\usepackage{graphicx}
\usepackage{epsfig}

\newcommand{\R}{\mathbb{R}}

\newcommand{\HL}{\dot{H}^1\times L^2}
\newcommand{\OR}{\overrightarrow}
\newcommand{\E}{\mathcal{E}}

\newcommand{\Fl}{{\rm Flux}}

\newcommand{\spartial}{ {\slash\!\!\! \partial} }
\newcommand{\ud}{u^{\dagger}}

\newcommand{\pua}{\partial^{\alpha}}

\newcommand{\wt}{\widetilde}
\numberwithin{equation}{section} 
\newcommand{\bt}{\overline{t}}
\newcommand{\br}{\overline{r}}
\newcommand{\bx}{\overline{x}}
\newcommand{\nablaxtu}{\left(\frac{|\nabla u|^2}{2}+\frac{|\partial_tu|^2}{2}\right)}
\newcommand{\flux}{\left(\frac{|\nabla u|^2}{2}+\frac{|\partial_tu|^2}{2}+\frac{x}{t}\cdot\nabla \ud\,\partial_tu\right)}
\newcommand{\fluxOne}{\left(\frac{|\nabla u|^2}{2}+\frac{|\partial_tu|^2}{2}+x\cdot\nabla \ud\,\partial_tu\right)}
\newcommand{\tr}{\tilde{r}}
\newcommand{\tu}{\widetilde{u}}
\newcommand{\bu}{\overline{u}}
\newcommand{\initial}{\|(u_0,\,u_1)\|_{\HL}}

\begin{document}
\title[]{Universality of blow up profile for small blow up solutions to the energy critical wave map equation}
\author{Thomas Duyckaerts$^1$}
\thanks{$^1$LAGA, Universit\'{e} Paris 13, Sorbonne Paris Cit\'{e}, UMR 7539 and Institut Universitaire de France. Partially supported by ERC Grant Blowdisol 291214 and French ANR Grant SchEq ANR-12-JS01-0005-01.}

\author{Hao Jia$^{2}$}
\thanks{$^2$IAS, Partially supported by NSF grant DMS-1600779, and grant DMS-1128155 through IAS.}

\author{Carlos Kenig$^{3}$}
\thanks{$^{3}$University of Chicago, Partially supported by NSF Grants DMS-1265429 and DMS-1463746.}

\author{Frank Merle$^4$}
\thanks{$^4$Cergy-Pontoise (UMR 8088), IHES. Supported by ERC Grant Blowdisol 291214.}%

\maketitle

{\bf Abstract.}  
In this paper we introduce the channel of energy argument to the study of energy critical wave maps into the sphere. More precisely, we prove a channel of energy type inequality for small energy wave maps similar to that in \cite{DJKM}, and as an application we show that for a wave map that has energy just above the degree one harmonic maps and that blows up in finite time,  the solution asymptotically de-couples into a regular part plus a traveling wave with small momentum, in the energy space. In particular, the only possible form of energy concentration is through the concentration of traveling waves. This is often called ``quantization of energy" at blow up. We also give a brief review of important background results in the subcritical and critical regularity theory for the two dimensional wave maps from \cite{Klainerman1,Klainerman2,TataruSterbenz,Tataru4,TaoSmallEnergy,KriegerSchlag}.

\begin{section}{Introduction}
We consider the Cauchy problem for wave map $u$ from $R^{2+1}$ with Minkowski metric to the standard 2-sphere $S^2\subset R^3$:
\begin{equation}\label{wavemap}
\partial_{tt}u-\Delta u=(|\nabla u|^2-|\partial_tu|^2)u, \,\,{\rm in}\,\,R^2\times [0,\infty),
\end{equation}
with initial data $\OR{u}(0):=(u,\partial_tu)(0)=(u_0,\,u_1)$. The initial data $(u_0,\,u_1)$ must satisfy the ``compatibility condition" that $|u_0|\equiv 1$ and $u_0\cdot u_1\equiv 0$. We shall only consider initial data that satisfies the compatibility condition. For simplicity, we assume that the initial data $\OR{u}(0)$ is smooth, $u_1$ is compactly supported, and that $u$ equals a fixed constant $u_{\infty}$ for large $x$. We call such wave maps classical, following the usual convention. Wave maps from the Minkowski space to a general Riemannian manifold $\mathcal{M}$ arise naturally as the hyperbolic counterpart of harmonic maps, and are given as critical points of the Lagrangian
$$\mathcal{L}(u):=\int_{R^3}|\nabla u|^2-|\partial_tu|^2\,dxdt,$$
for $u:R^3\to \mathcal{M}.$ It is sometimes more convenient to adopt the more geometric notation: set for $\alpha=0,1,2$ that $\partial_{\alpha}=\partial_t$ if $\alpha=0$, $\partial_{\alpha}=\partial_{x_{\alpha}}$ if $\alpha=1,\,2$, and that $\partial^{\alpha}=-\partial_{\alpha}$ if $\alpha=0$, $\partial^{\alpha}=\partial_{\alpha}$ if $\alpha=1,\,2$. This is of course just using the Minkowski metric to lower or upper the index. We adopt the Einstein summation convention with repeated indices and view $u$ as a column vector. We also use the standard notation that $x^0=x_0=t,\,x^j=x_j$ for $j=1,\,2$. Then equation (\ref{wavemap}) can be written as
\begin{equation}\label{Mainwavemap}
-\,\partial_{\alpha}\partial^{\alpha}u=u\,\partial^{\alpha}\ud\partial_{\alpha}u,
\end{equation}
where $\ud$ is the transpose of $u$.\\

The wave map equation has been intensively studied, as a natural geometric wave equation and as models from physics - including general relativity and gauge theories. The study of the Cauchy problem and the dynamics of solutions in the equivariant setting was initiated in the works of Shatah and Tahvildar-Zadeh \cite{Shatah1,Shatah2}, Christodoulou and Tahvildar-Zadeh \cite{Chri,ChriTah}, and Struwe \cite{Struwe2}. In general, the wave map can also develop a singularity in finite time by concentrating energy in a small region.  Indeed, singular solutions in the form of a shrinking soliton plus a residue term have been constructed for the $2+1$ equivariant dimensional wave map equation by Krieger, Schlag and Tataru \cite{KSTwavemap} with prescribed rate, by Rodnianski and Sterbenz \cite{RodSter} in a stable regime for high equivariance wave maps, and by Raphael and Rodnianski \cite{RapRod} for co-rotational wave maps. We also refer to the recent survey \cite{SchlagICM} for further discussion. The Cauchy 
problem for the wave maps without equivariant symmetry is more complicated.  Recall that equation (\ref{wavemap}) is invariant under the natural scaling
\begin{equation}\label{eq:scaling}
u\to u_{\lambda}(x,t)=u(\lambda x,\,\lambda t),\,\,\,\,(u_0,\,u_1)\to (u_{0\lambda }(x),\,u_{1\lambda }(x))=(u_0(\lambda x),\,\lambda u_1(\lambda x)),
\end{equation}
and the conserved energy
\begin{equation}\label{energy}
\mathcal{E}(\OR{u}):=\int_{R^2}\nablaxtu(x,t)\,dx
\end{equation}
is invariant under the scaling (\ref{eq:scaling}). Scale-invariance of the equation plays an essential role in both the Cauchy problem and the dynamics of solutions. We note that for equation (\ref{wavemap}), the natural initial data space invariant under the scaling (\ref{eq:scaling}) is the energy space $\HL$, and hence the equation is called {\it energy critical}.
The works of Klainerman and Machedon \cite{KlainermanMachedon1,KlainermanMachedon2,KlainermanMachedon3}, and subsequently Klainerman and Selberg \cite{Klainerman1,Klainerman2}, and Selberg \cite{Selberg} established wellposedness in the subcritical space $\dot{H}^{s-1}\times H^{s-1}$ with $s>1$, and introduced important ideas on the bilinear and null form estimates that also played an important role in the critical theory. The Cauchy problem for the wave map equation in the critical space $\HL$ is more difficult, and was addressed  in the breakthrough work of Tao \cite{TaoSmallEnergy}, using the null frame spaces introduced by Tataru \cite{Tataru1} and Tao's idea of gauge transform \cite{TaoHighDimension}. The global wellposedness for the energy critical wave maps has also been intensively studied, see the works of Tao \cite{taoIII,TaoIV,TaoV,TaoVI}, Sterbenz and Tataru \cite{TataruSterbenz,Tataru4}, and Krieger and Schlag \cite{KriegerSchlag}. In particular  \cite{TataruSterbenz,Tataru4} proves that if a 
wave map blows up in finite time, then after suitable transformation using symmetry, it must converge locally to a harmonic map.  Recently, Grinis \cite{Grinis} extended the result of \cite{TataruSterbenz,Tataru4} and obtained a more complete characterization of the wave map {\it strictly} inside the lightcone along a sequence of times. He showed that all the energy concentration along a sequence of times strictly inside the lightcone must be in the form of traveling waves, by showing that there is no energy in the so called ``neck region". His work will also be important for our purposes. There are also many works on the study of large equivariant wave maps in connection to the ``soliton resolution conjecture", see for example \cite{Cotemap1,Cotemap2,KLS,KLIS1,KLIS2}. The new ingredient in these works is the {\it channel of energy inequality} first introduced in \cite{DKM,DKMsmall}, which provides strong decoupling mechanism between the dispersion and solitary waves. The channel of energy type 
inequality 
has been applied to many other problems, see for example \cite{Kenig4dWave,DKMsuper,Casey,JiaLiuXu,JLSX} for application to semilinear energy critical wave equations. These channel of energy inequalities in many cases depend crucially on the radial assumption and are sensitive to the dimensions. \\

Despite significant progresses in the study of energy critical wave maps, there are still many interesting and deep questions that remain. For example, in the equivariant setting ( in \cite{Cotesoliton} for co-rotational and \cite{JiaKenig} for all equivariant wave maps) it was proved that along a sequence of times the solution asymptotically de-couples into a finite sum of harmonic maps plus a regular part in the finite time blow up case or a linear wave in the global existence case.  A natural question is then if such a decomposition persists for all times. ({\it Soliton resolution conjecture} predicts that (for generic targets) this decomposition persists for all times.) If one imposes certain energy constraint that effectively rules out multi-soliton configuration, then the answer is yes. But in general the question remains open. To answer this question, it seems that one needs to understand the interaction of solitons that are separated either by scales, or by distances, which appears to be a 
challenging task. One can also ask what happens if we remove the equivariance assumption. A natural step seems to be the extension of the corresponding result in the equivariant case, that is, to prove the soliton resolution conjecture along a sequence of times. In this case, a new difficulty appears that is not present in the equivariant setting. In the equivariant setting, it was known (see \cite{Chri,ChriTah}) that there is asymptotically no energy accumulation in the so called ``self similar" region. In particular, in the equivariant case there can not be any energy concentration near the boundary of the singularity lightcone $|x|<T_+-t$ as $t\to T_+$, assuming that the solution blows up at time $T_+$. As far as the authors know, it is an open question how to rule out energy concentration near the boundary of the singularity lightcone in the general case. \\

Our work addresses the question of ruling out energy concentration near the boundary of the singularity lightcone, in the restricted case where the energy is only slightly higher than the energy of the degree one co-rotational harmonic maps. We believe that the methods used here apply to wave maps into more general targets without any size restrictions. However it appears that one has to overcome serious obstacles in the current perturbative setup to achieve this goal. We also believe that the case $T_+=\infty$, and the solution does not scatter, can also be addressed by the methods developed here. We plan to address these questions in future work.\\

Let us briefly summarize our main results.\\

Our first main goal is to introduce the channel of energy argument to the study of wave map equations. The  channel of energy inequality for outgoing waves (see Theorem \ref{th:channelofenergywavemapmain}) has played an essential role in the recent proof of soliton resolution conjecture along a sequence of times in \cite{DJKM}. Unlike previously known channel of energy inequalities, the version for outgoing waves turns out to be rather robust in that it works for nonradial solutions in all dimensions. The outgoing condition is natural. For instance, any linear wave at large time will satisfy such outgoing conditions. More interestingly, in the blow up case and away from the concentrating solitons, the dispersed energy that might concentrate near the boundary of the singularity lightcone also satisfies the outgoing condition. Thus this channel of energy inequality appears to be applicable to a wide range of problems. Of course, in comparison with the focusing energy critical wave equation considered in \cite{
DJKM}, the wave map equations  are much more complicated and the current perturbation results are not as precise as in the case of focusing energy critical  wave equations. At this time we can only extend the results from \cite{DJKM} partially, and prove the channel of energy inequality for small data. 
\begin{theorem}\label{th:channelofenergywavemapmain}
Fix $\beta\in (0,1)$. There exists a small  $\delta=\delta(\beta)>0$ and sufficiently small  $\epsilon_0=\epsilon_0(\beta)>0$,  such that if $u$ is a classical wave map with energy $\mathcal{E}(\OR{u})<\epsilon_0$ satisfying
\begin{equation}\label{eq:outgoingconditionmain}
\|(u_0,\,u_1)\|_{\HL\left(B^c_{1+\delta}\cup B_{1-\delta}\right)}+\|\spartial u_0\|_{L^2}+\|\partial_ru_0+u_1\|_{L^2}\leq \delta \|(u_0,\,u_1)\|_{\HL},
\end{equation}
then for all $t\ge 0$, we have
\begin{equation}\label{eq:smallchannelMain}
\int_{|x|>\beta+t}|\nabla_{x,t}u|^2(x,t)\,dx\ge \beta\,\|(u_0,\,u_1)\|_{\HL}^2.
\end{equation} 
\end{theorem}

As an application for the channel of energy inequality (\ref{eq:smallchannelMain}), we obtain the following classification of finite time blow up wave maps $u$ with energy 
\begin{equation}\label{eq:energyconstraint}
\E(\OR{u})<\E(Q,0)+\epsilon_0^2,
\end{equation} where $Q$ is the harmonic map with the least energy (which is equal to $4\pi$). Denote $\mathcal{M}_1$ as the space of degree one harmonic maps, \footnote{These harmonic maps are all co-rotational with respect to certain axis of symmetry.} and let
$$\mathcal{M}_{\ell,1}:=\left\{Q_{\ell}:\,Q\in \mathcal{M}\right\},$$
where
\begin{equation}\label{eq:QlI}
Q_{\ell}(x,t)=Q\left(x-\frac{\ell\cdot x}{|\ell|^2}\ell+\frac{\frac{\ell\cdot x}{|\ell|^2}\ell-\ell t}{\sqrt{1-|\ell|^2}}\right),
\end{equation}
is the Lorentz transformation of the harmonic map $Q$. Then we have

\begin{theorem}\label{th:smallsolitonresolutionMain}
Let $u$ be a classical wave map with energy $\E(\OR{u})<\E(Q,0)+\epsilon_0^2$, that blows up at a finite time  $T_+$ and at the origin. Assume that $\epsilon_0$ is sufficiently small. Then there exists $\ell\in R^2$ with $|\ell|\ll 1$, $x(t)\in R^2,\,\lambda(t)>0$ with 
$$\lim_{t\to T_+}\frac{x(t)}{T_+-t}=\ell,\,\,\,\lambda(t)=o\left(T_+-t\right),$$
and $(v_0,\,v_1)\in \HL\cap C^{\infty}(R^2\backslash\{0\})$ with $(v_0-u_{\infty},\,v_1)$ being compactly supported, such that
\begin{eqnarray*}
{\rm (i)}&&\inf\bigg\{\left\|\OR{u}(t)-(v_0,\,v_1)-\left(Q_{\ell},\,\partial_tQ_{\ell}\right)\right\|_{\HL}:\,Q_{\ell}\in \mathcal{M}_{\ell,1}\bigg\}\to 0,\,\,{\rm as}\,\,t\to T_+;\\
&&\\
{\rm (ii)}&&\bigg\|\OR{u}(t)-(v_0,\,v_1)\bigg\|_{\HL\left(R^2\backslash B_{\lambda(t)}(x(t))\right)}\to 0\,\,{\rm as}\,\,t\to T_+.
\end{eqnarray*}
\end{theorem}

Heuristically speaking, the above theorem says that at blow up time, the wave map essentially consists of two parts, one regular part outside the lightcone $|x|>T_+-t$, and a traveling wave with small velocity $\ell$ that concentrates in a small region (in comparison with the size of the cone) near the point $\ell (T_+-t)$. In addition, there are no other types of energy concentration. It is an interesting question to ask about the {\it finer} dynamics of the traveling wave in the region $|x-x(t)|<\lambda(t)$, such as how the axis of rotation of the wave map evolves. It is not clear to us at this moment if the axis of rotation could fail to stabilize. We believe though that it is impossible, at the level of energy regularity, to say more about the finer dynamics of the scale $\lambda(t)$ and the center $x(t)$, due to the symmetries of the equation.\\

Let us very briefly explain the strategy of the proof. The proof of the channel of energy inequality in Theorem \ref{th:channelofenergywavemapmain} uses the extension of the linear channel of energy inequality for outgoing waves from \cite{DJKM} to two dimensions. For the wave maps however, we need to also show that these outgoing conditions are in some sense stable (for {\it most frequencies}) with respect to frequency projections, in order to use the perturbative results for wave maps which deal with each frequency piece of the map separately. To prove Theorem \ref{th:smallsolitonresolutionMain}, let us take a wave map as in Theorem \ref{th:smallsolitonresolutionMain}. Then by the result of Tataru and Sterbenz \cite{Tataru4}, along a sequence of times, we can extract a traveling wave from the wave map. A little more effort also shows that there is no other possible energy concentration strictly inside the lightcone except in the neck region, thanks to the energy constraint (\ref{eq:energyconstraint}). By 
Grinis' result \cite{Grinis} there is no energy in the neck region either. Thus all residue energy has to concentrate near the boundary of the singularity lightcone $|x|<T_+-t$. In addition, such residue energy has to be small, again thanks to the energy constraint. We apply the channel of energy inequality to rule out this residue energy. This is a crucial step and the main new point of our paper. Hence inside the lightcone (not only {\it strictly inside} the lightcone) the amount of energy is asymptotically just the energy of the traveling wave. Then by the  coercivity of energy near the traveling wave, we conclude that in fact the wave map is trapped in smaller and smaller neighborhoods of the traveling wave, and thus has to stay close to the traveling wave for all times $t<T_+$, not just along a sequence of times. This completes the proof of the main Theorem \ref{th:smallsolitonresolutionMain}.\\

Our paper is organized as follows:
\begin{itemize}
\item In section 2, we recall the necessary subcritical and critical regularity results for the wave equation;

\item In section 3, we prove the channel of energy inequality for small wave maps;

\item In section 4, we recall the Morawetz estimates;

\item In section 5, we prove the decomposition into regular part and traveling wave along a sequence of times;

\item In section 6, we prove certain coercive property of energy in the neighborhood of the traveling wave, and establish the decomposition for all times;


\end{itemize}
Throughout the paper, we shall use the notation
$$\left\|f\right\|_{\dot{H}^1(E)}:=\left\|\nabla f\right\|_{L^2(E)},$$
for any measurable set $E$. If $s>1$, we will write
$$\left\|f\right\|_{\dot{H}^s}:=\left\||D|^sf\right\|_{L^2(\R^2)},$$
where $|D|^s$ is the Fourier multiplier with symbol $|\xi|^s$, and say that a distribution $f$ is in $\dot{H}^s$ when $f\in H^s_{\rm{loc}}(\R^2)$ and the above seminorm is finite.
\end{section}

\begin{section}{Preliminaries}
In this section, we briefly review the subcritical and critical regularity results for the two dimensional wave maps into the sphere, that will be needed below. 
\subsection{Local wellposedness in $H^s$ for $s>1$}
It is well known from the works of Klainerman and Machedon \cite{KlainermanMachedon1,KlainermanMachedon2,KlainermanMachedon3}, Klainerman and Selberg \cite{Klainerman1,Klainerman2}, and Selberg \cite{Selberg}, that the wave map equation (\ref{Mainwavemap}) is locally wellposed in $H^s\times H^{s-1}$ for $s>1$. In this subsection we recall the necessary regularity results from these works without giving proofs and refer the reader to the above cited works, and especially the survey \cite{Klainerman2} for details.\\

 Since the spaces in which one can prove existence and uniqueness involves spacetime Fourier transforms even when one only considers local in time solutions, we have to be more precise on the Banach spaces which are used to hold the solutions and the nonlinearities. \\
We shall denote $\mathcal{F}(u)$ as the spacetime Fourier transform of $u$.
For $s,\,b\in R$, and tempered distribution $u\in \mathcal{S'}(R^3)$, define
\begin{equation}\label{HyperbolicSobolev}
\|u\|_{X^{s,b}(R^3)}:=\left(\,\int_{R^3}(1+|\xi|^2)^{s}\,(1+||\xi|-|\tau||)^{2b}\,|\mathcal{F}(u)(\xi,\tau)|^2\,d\xi d\tau\,\right)^{\frac{1}{2}},
\end{equation}
and set
\begin{equation}\label{HyperbolicSpace}
X^{s,b}(R^3):=\left\{u\in\mathcal{S'}(R^3):\,\|u\|_{X^{s,b}(R^3)}<\infty\right\}.
\end{equation}

We record the following wellposedness result for equation (\ref{Mainwavemap}) in the subcritical space $\dot{H}^s\times H^{s-1}$. We shall always assume that the initial data $(u_0,\,u_1)$ for (\ref{Mainwavemap}) satisfies the ``admissibility condition" that $|u_0|\equiv 1$ and $\ud_0\,u_1\equiv 0$.
\begin{theorem}\label{th:subcriticallocalwellposedness}
For $s>1$ and $\frac{1}{2}<b<\min\{s-\frac{1}{2},\,1\}$. Suppose that $(u_0,\,u_1)\in \dot{H}^{s}\times H^{s-1}$ and that $u_0$ equals a constant $u_{\infty}\in S^2$ for large $x$. Then for $T=T(\|(u_0-u_{\infty},\,u_1)\|_{H^s\times H^{s-1}})>0$ sufficiently small, there exists a unique solution $u$ to equation (\ref{Mainwavemap}) with initial data $(u_0,\,u_1)$ on $R^2\times (-T,T)$ in the sense of distributions, which satisfies the following properties
\begin{eqnarray*}
&(1)& u-u_{\infty}\in C(I, H^s\times H^{s-1});\label{eq:wellpro1}\\
&(2)& {\rm there\,\,exists}\,\,\overline{u}\in L^2(R^3)\,\,{\rm with}\,\,\overline{u}|_{R^2\times I}\equiv u-u_{\infty} \,\,{\rm and\,\,}\nabla_{x,t}\overline{u}\in X^{s-1,b},\label{eq:wellpro2}
\end{eqnarray*}
where $I=(-T,\,T).$
\end{theorem}

\smallskip
\noindent
{\it Remark.} The above theorem provides a rigorous definition of solutions to equation (\ref{Mainwavemap}). (2) is important, as (1) by itself is not sufficient to guarantee uniqueness when $s$ is close to $1$. One could of course choose to work directly with smooth wave maps, instead of these low-regularity wave maps. However, below we shall need to extend a locally (in space) defined map to a global one, and it is much more convenient to have such extensions in the framework of $H^s$ solutions, rather than smooth solutions.\\

Solutions from Theorem \ref{th:subcriticallocalwellposedness} can be extended to a maximal interval of existence, more precisely, we have
\begin{corollary}\label{cor:maximalextension}
For $s>1$. Suppose that $(u_0,\,u_1)\in \dot{H}^{s}\times H^{s-1}$ and that $u_0$ equals a constant $u_{\infty}\in S^2$ for large $x$. Then there exists $T_+\in(0,\infty],\,T_{-}\in [-\infty,0)$, such that for any $T_{-}<T_1<T_2<T_{+}$, $u$ is a distributional solution to equation (\ref{Mainwavemap}), satisfying (1) and (2) on $I=(T_1,\,T_2)$, and that if $T_{+}<\infty$, then
\begin{equation}\label{eq:blowupcriteria}
\lim_{t\to T_+}\|\OR{u}(t)\|_{\dot{H}^s\times H^{s-1}}=\infty.
\end{equation}
Similar conclusion holds for $T_{-}$. Such $u$ is unique. In addition, if $(u_0,\,u_1)\in \dot{H}^{s_1}\times H^{s_1-1}$ for some $s_1>s$, then $u$ satisfies (\ref{eq:wellpro1}) and (\ref{eq:wellpro2}) with $s$ being replaced by $s_1$ on any $I=(T_1,\,T_2)\Subset (T_{-},\,T_{+})$. $T_+$ and $T_{-}$ are called the maximal time of existence for the solution $u$. 
\end{corollary}

\subsection{Critical wellposedness results} Perhaps not surprisingly, our work depends crucially on the regularity results of Tao \cite{TaoHighDimension,TaoSmallEnergy}, Tataru \cite{Tataru1}, and Sterbenz-Tataru \cite{Tataru4,TataruSterbenz}. See also the work of Krieger and Schlag \cite{KriegerSchlag}. In this section we recall some important results for wave maps in the energy space from \cite{Tataru1,TaoSmallEnergy,TataruSterbenz}, that will be needed below.\\

 In order to control the solution at the $\dot{H}^1\times L^2$ level of regularity, we need to use more sophisticated spaces. The precise definitions of these spaces are not very important for us, but we shall need the following properties that we briefly review below. \\

Fix a radial function $\Phi\in C_c^{\infty}(R^2)$ with $\Phi|_{B_1}\equiv 1$ and ${\rm supp}\,\Phi\Subset B_2$. Let $\Psi(x):=\Phi(x)-\Phi(2x)$, and $\Psi_k(x)=\Psi(x/2^k)$ for each $k\in \mathbb{Z}$. Then ${\rm supp}\,\Psi\Subset B_2\backslash B_{1/2}$, and
$$\sum_{k\in \mathbb{Z}}\Psi_k\equiv 1,\,\,\,\,{\rm for}\,\,|\xi|\neq 0.$$
Recall that the Littlewood-Paley projection $P_k$ and $P_{<k}$ are defined as
$$\widehat{P_kf}(\xi)=\Psi_k(\xi)\,\widehat{f}(\xi),$$
and 
$$P_{<k}f=\sum_{k'<k}P_{k'}f.$$
We will also use the notations $u_k:=P_ku$ and $u_{<k}=P_{<k}u$.
Then
$$\sum_{k\in\mathbb{Z}}P_kf=f,$$
for all $f\in L^2(R^2)$. We use the same definitions as in \cite{TaoSmallEnergy} for the spaces $S[k],\,N[k]$, which are translation invariant Banach spaces of distributions on $\R^2_x\times \R_t$ containing Schwartz functions whose partial Fourier transform in the $x$ variable is supported in $\{2^{k-3}\leq |\xi|\leq 2^{k+3}\}$, $\{2^{k-4}\leq |\xi|\leq 2^{k+4}\}$ respectively. For each $k$, we shall use the space $S[k]$ to hold the frequency localized piece $P_ku$ of the solution $u$, and use the space $N[k]$ to hold the frequency localized piece $P_kf$ of the nonlinearity $f:=u\,\partial_{\alpha}\ud\partial^{\alpha}u$. Define the $S(1)$ norm as
\begin{equation}
\|f\|_{S(1)}:=\|f\|_{L^{\infty}}+\sup_k\|P_kf\|_{S[k]}.
\end{equation}
The spaces $S[k]$ and $N[k]$ satisfy the following properties.
\begin{theorem}\label{th:propertySN}
There exists a small universal constant $\kappa>0$, such that\\
(1)\,{\bf (Algebra property)}\, For Schwartz functions $\phi,\,\psi$ with $\psi\in S[k_2]$, we have
\begin{equation}\label{eq:algebra}
\|P_k(\phi\psi)\|_{S[k]}\lesssim 2^{-\kappa (k_2-k)_{+}}\|\phi\|_{S(1)}\|\psi\|_{S[k_2]};
\end{equation}
(2)\,{\bf (Product property)}\, For Schwartz functions $f,\,\psi$ with $f\in N[k_2]$, we have
\begin{equation}\label{eq:product}
\|P_k(f\psi)\|_{N[k]}\lesssim 2^{-\kappa (k_2-k)_{+}}\|\psi\|_{S(1)}\|f\|_{N[k_2]};
\end{equation}
(3)\,{\bf (Null form estimate)} \,For Schwartz functions $\phi,\,\psi$ with $\phi\in S[k_1]$, $\psi\in S[k_2]$, we have
\begin{equation}\label{eq:null}
\|P_k(\partial^{\alpha}\phi\,\partial_{\alpha}\psi)\|_{N[k]}\lesssim 2^{-\kappa (\max\{ k_1,\,k_2\}-k)_{+}}\|\phi\|_{S[k_1]}\|\psi\|_{S[k_2]};
\end{equation}
(4)\,{\bf (Trilinear estimate)} \,For Schwartz functions $\phi,\,\varphi,\,\psi$ with $\phi\in S[k_1]$, $\varphi\in S[k_2]$ and $\psi\in S[k_3]$, we have
\begin{multline}\label{eq:trilinear}
\|P_k(\phi\,\partial^{\alpha}\varphi\partial_{\alpha}\psi)\|_{N[k]}\\
\lesssim 2^{-\kappa (k_1-\min\{k_2,\,k_3\})_+}2^{-\kappa(\max\{k_1,k_2,k_3\}-k)_+}\|\phi\|_{S[k_1]}\|\varphi\|_{S[k_2]}\|\psi\|_{S[k_3]}.
\end{multline}
(5)\,{\bf (Linear wave estimate)} \,For solution $u^L$ to the linear wave equation $$\partial_{tt}u^L-\Delta u^L=f,$$ with initial data $(u_0,\,u_1)$, we have
\begin{equation}\label{eq:inhomogeneousenergyestimates}
\|P_ku^L\|_{S[k]}\lesssim \|P_k(u_0,\,u_1)\|_{\HL}+\|P_kf\|_{N[k]}.
\end{equation}
(6)\,{\bf ($S[k]$ controls energy)} \,For $u\in S[k]$, we have
\begin{equation}\label{eq:Senergy}
\|\nabla_{x,t}u\|_{L^{\infty}_tL^2_x}\lesssim \|u\|_{S[k]}.
\end{equation}

\end{theorem}

\smallskip
\noindent
{\it Remark.} These estimates were proved in \cite{TaoSmallEnergy}, and some of them are slightly more general than those stated in the main summary of the properties of $S[k],\,N[k]$ in Theorem 3  from \cite{TaoSmallEnergy}. However, they can be found elsewhere in that paper. More precisely, the algebra estimate (\ref{eq:algebra}) is  a consequence of equation (125) and (126) at page 516; the product estimate (\ref{eq:product}) is a consequence of (119) at page 510;  the null form estimate (\ref{eq:null}) is (134) at page 523; the trilinear estimate (\ref{eq:trilinear}) is taken from the first formula at page 529. We also note that $P_{k'}$ is bounded from $S[k]$ to $S[k]$, by the translation invariance of the Banach space $S[k]$. (6) implies that $\|u\|_{L^{\infty}}\lesssim \|u\|_{S[k]}$. Another useful property of $S[k]$ is the weak stability of $S[k]$: if $u_i\to u$ in the sense of distributions and $u_i\in S[k]$ with $\|u_i\|_{S[k]}\leq 1$, then $\|u\|_{S[k]}\leq 1$. See a similar statement in (vii) of 
page 323 in \cite{TataruRough}. We will use these estimates extensively below.\\

Tao \cite{TaoSmallEnergy} introduced a very useful notation to keep track of multilinear expressions. More precisely, for scalar functions $\phi_1,\dots,\phi_l$, we use $L\left(\phi_1,\dots,\phi_l\right)$ to denote multilinear expression of the form
\begin{equation*}
L\left(\phi_1,\dots,\phi_l\right):=\int K\left(y_1,\dots,y_l\right)\phi_1(x-y_1)\cdots\phi_l(x-y_l)\,dy_1\cdots dy_l,
\end{equation*}
with a measure $K$ of bounded mass. In many cases, $\phi_1,\dots,\phi_l$ could also be expressions involving components $\phi^{j_1}_1,\dots,\phi^{j_l}_l$ and in such cases, we also assume that $K$ depends on $j_1,\dots,j_l$, but for the ease of notations, we shall suppress this dependence. By the translation invariance of the spaces $S[k],\,N[k]$, the estimates in Theorem \ref{th:propertySN} extend to expressions of the form $L(\phi,\psi),\,L(\partial^{\alpha}\phi,\partial_{\alpha}\psi)$ instead of just $\phi\,\psi$ and $\partial^{\alpha}\phi\,\partial_{\alpha}\psi$.\\

 Let us record here the following useful Lemma from \cite{TaoSmallEnergy}.
\begin{lemma}\label{lm:commutinglemma}
For Schwartz functions $f,\,g$, we have
\begin{equation}
P_k(f\,g)-P_kf\cdot g=2^{-k}L\left(f,\,\nabla g\right).
\end{equation}
\end{lemma}

\smallskip
\noindent
{\it Proof.} This Lemma is taken from \cite{TaoSmallEnergy}, we include the short proof for the convenience of readers. We have
\begin{eqnarray*}
&&P_k(f\,g)(x)-P_kf(x)\,g(x)\\
&&=\int 4^k\check{\Psi}(2^ky)f(x-y)g(x-y)\,dy-\int 4^k\check{\Psi}(2^ky)f(x-y)g(x)\,dy\\
&&=\int_0^1\int_{R^2}-4^k\check{\Psi}(2^ky)f(x-y)\,y_j\partial^jg(x-ty)\,dtdy\\
&&=-2^{-k}\int_0^1\int_{R^2}4^k(2^ky_j)\,\check{\Psi}(2^ky)f(x-y)\,\partial^jg(x-ty)\,dtdy\\
&&=2^{-k}L(f,\,\nabla g).
\end{eqnarray*}
The proof is complete.\\

Let us recall the definition of {\it frequency envelop} introduced in \cite{TaoSmallEnergy}. Fix positive $\vartheta$ such that $\vartheta\leq \frac{\kappa}{100}$, where $\kappa$ is as in Theorem \ref{th:propertySN}.

\begin{definition}
$(c_k)\in \ell^2$ is called a frequency envelop if $c_k>0$ and $c_{k_1}\leq 2^{\vartheta |k_1-k_2|}c_{k_2}$.
\end{definition}
For any frequency envelop $c=(c_k)$, define the norm $S(c)$ as
\begin{equation}
\|\phi\|_{S(c)}:=\|\phi\|_{L^{\infty}}+\sup_kc_k^{-1}\|P_k\phi\|_{S[k]},
\end{equation}
and the space $S(c)$ as 
\begin{equation}
S(c):=\{f\in L^{\infty}:\,\|f\|_{S(c)}<\infty\}.
\end{equation}
 Note that $1\in S(c)$. The main property of the space $S(c)$ that we shall use below is that $S(c)$ a Banach algebra. 
\begin{lemma}
$S(c)$ is a Banach algebra.
\end{lemma}

\smallskip
\noindent
{\it Proof:} This was proved in \cite{TaoSmallEnergy}. We include the short proof for the convenience of readers. We need to prove 
\begin{equation}
\|\phi\,\psi\|_{S(c)}\lesssim \|\phi\|_{S(c)}\,\|\psi\|_{S(c)}.
\end{equation}
We note that $\|\phi\|_{S(1)}\lesssim_c\|\phi\|_{S(c)}$ and $\|\phi_{<k}\|_{S(c)}\lesssim \|\phi\|_{S(c)}$. We can normalize $\|\phi\|_{S(c)}=\|\psi\|_{S(c)}=1$.
For each $k\in \mathbb{Z}$, we have 
\begin{eqnarray*}
&&\left\|P_k(\phi\,\psi)\right\|_{S[k]}\\
&&=\left\|P_k\left(\phi_{>k-10}\psi\right)+P_k\left(\phi_{\leq k-10}\,\psi_{>k-10}\right)+P_k\left(\phi_{\leq k-10}\psi_{\leq k-10}\right)\right\|_{S[k]}.
\end{eqnarray*}
Note that $$P_k\left(\phi_{\leq k-10}\psi_{\leq k-10}\right)\equiv 0.$$
We get that
\begin{eqnarray*}
\left\|P_k(\phi\,\psi)\right\|_{S[k]}&\lesssim& \sum_{k_1>k-10}\left\|P_k\left(P_{k_1}\phi\,\psi\right)\right\|_{S[k]}+\sum_{k_2>k-10}\left\|P_k\left(\phi_{\leq k-10}\,P_{k_2}\psi\right)\right\|_{S[k]}\\
&\lesssim&\sum_{k_1>k-10}2^{-\kappa (k_1-k)_+}\left\|P_{k_1}\phi\right\|_{S[k_1]}\,\|\psi\|_{S(1)}\\
&&\hspace{.3in}+\,\sum_{k_2>k-10}2^{-\kappa (k_2-k)_+}\left\|P_{k_2}\psi\right\|_{S[k_2]}\,\|\phi_{\leq k-10}\|_{S(1)}\\
&\lesssim&\sum_{k_1>k-10}2^{-\kappa (k_1-k)}c_{k_1}+\sum_{k_2>k-10}2^{-\kappa (k_2-k)}c_{k_2}\\
&\lesssim&\sum_{k'>k-10}2^{-(\kappa-\vartheta)(k'-k)}c_{k}\lesssim c_k,
\end{eqnarray*}
and this finishes the proof.\\

Let us recall the following global wellposedness theorem for wave maps from Tao \cite{TaoSmallEnergy}. 
\begin{theorem}\label{th:globalregularity}
There exists an $\varepsilon>0$ sufficiently small such that the following is true. Suppose that $(u_0,u_1)$ is smooth, $u_0-u_{\infty},\,u_1$ are compactly supported, and that $u_1^{\dagger}\cdot u_0\equiv 0$. Assume that $\left(\|P_k(u_0,\,u_1)\|_{\HL}\right)$ lies under a frequency envelop $c=(c_k)$ \footnote{For non-negative sequences $(a_k)$ and $(b_k)$, we say that $(a_k)$ lies below $(b_k)$ if $a_k\leq b_k$ for each $k$.} with 
$$\|c_k\|_{\ell^2}\leq \varepsilon.$$
Then the wave map $u$ with initial data $(u_0,u_1)$ is global, and moreover
\begin{equation}\label{eq:GlobalControl}
\|P_ku\|_{S[k]}+\sup_{t\in R}\|P_k\OR{u}(t)\|_{\HL}\leq Cc_k,
\end{equation}
for some universal $C$.
\end{theorem}

\smallskip
\noindent
{\it Remark:} By approximations by smooth maps, and the wellposedness for equation (\ref{Mainwavemap}) in $\dot{H}^s\times H^{s-1}$ for $s>1$, we can relax the smoothness requirement for the initial data in the above theorem to $(u_0,\,u_1)\in \dot{H}^s\times H^{s-1}$.\\

Fix $\epsilon_{\ast}>0$ be sufficiently small, so that classical wave maps with energy smaller than $C\epsilon_{\ast}$ exists globally for a sufficiently large universal $C>1$.
In later sections, we shall need the following local-in-space smoothness result, when the initial data is locally but not globally smooth.
\begin{lemma}\label{lm:localsmoothness}
Let $(u_0,\,u_1)\in \dot{H}^s\times H^{s-1}$ for some $s>1$ and that $(u_{0}-u_{\infty},\,u_1)$ is compactly supported, with $|u_0|\equiv 1$ and $u_0\cdot u_1\equiv 0$. Assume that $(u_0,\,u_1)$ is smooth in $B_1$, and that $\initial\leq\epsilon_{\ast}.$ Then the global solution $u$ is smooth in $\{(x,t):\,|x|<1-|t|\}$.
\end{lemma}

\smallskip
\noindent
{\it Proof:} By the small energy global existence result and the subcritical Cauchy theory, we get that
\begin{equation}\label{eq:Hsnorm}
\sup_{0\leq t<1}\|u(t)\|_{\dot{H}^s}\leq C\left( \|(u_0-u_{\infty},\,u_1)\|_{H^s\times H^{s-1}}\right).
\end{equation}
Denote $\lambda:=C\left( \|(u_0-u_{\infty},\,u_1)\|_{H^s\times H^{s-1}}\right)$. Fix $\br>0$ small, and take $B_{\br}(\bx)\subset B_1$. Set
$$\bu_0=\frac{1}{\pi\br^2}\int_{B_{\br}(\bx)}u.$$
Then by Sobolev inequality, we obtain that
\begin{equation}\label{eq:smallinsmallball}
\left\|u_0-\bu_0\right\|_{L^{\infty}\left(B_{\br}(\bx)\right)}\lesssim_{\lambda}\br^{s-1}.
\end{equation}
Take a smooth cutoff function $\eta$ such that $\eta\equiv 1$ in $B_{\br-\br^{1+\delta}}(\bx)$ with some $\delta\in (0,\,2(s-1))$, and $\eta\equiv 0$ outside $B_{\br}(\bx)$. In addition, we can require that 
\begin{equation}\label{eq:smalleta}
|\nabla \eta|\lesssim \br^{-1-\delta}.
\end{equation}
Define
\begin{equation*}
\left(\tu_0,\,\tu_1\right)=\left(P\left[\eta\,(u_0-\bu_0)+\bu_0\right],\,\eta\, u_1\right),
\end{equation*}
where for each vector $v\neq 0$
$$Pv=\frac{v}{|v|}.$$
Then $\tu_0,\,\tu_1$ are smooth, and 
\begin{equation}\label{eq:localequal11}
\left(\tu_0,\,\tu_1\right)\equiv (u_0,\,u_1),\,\,\,\,{\rm in}\,\,B_{\br-\br^{1+\delta}}(\bx).
\end{equation}
Moreover, we can verify by direct computation thanks to (\ref{eq:smallinsmallball}) and (\ref{eq:smalleta}) that
$$\left\|\left(\tu_0,\,\tu_1\right)\right\|_{\HL}\lesssim\epsilon_{\ast},$$
if $\br$ is chosen sufficiently small.
Hence the solution $\tu$ to the wave map equation with the initial data $\left(\tu_0,\,\tu_1\right)$ is smooth and global. By (\ref{eq:localequal11}), $u\equiv \tu$ for $|x-\bx|<\br-\br^{1+\delta}-|t|$, and is thus smooth for $|x-\bx|<\br-\br^{1+\delta}-|t|$. By moving around $\bx$ and finite speed of propagation, we conclude that $u$ is smooth in $\{(x,t):\,|x|<1-2\br^{1+\delta}-|t|,\,|t|<\br\}$. We can apply the same technique at $|t|=\br,\,2\br$ and so on, and conclude recursively that $u$ is smooth in $\{(x,t):\,|x|-2k\br^{1+\delta}-|t|,\,|t|<k\br\}$ for $k=1,\,2,\dots$ with $(k+2)\br<1.$
Hence, $u$ is smooth in $\{(x,t):\,|x|<1-C\br^{\delta},\,\,|t|<1-3\br\}$. Since $\br$ can be taken arbitrarily small, the lemma follows.\\

Since the global regularity result for small energy requires that the initial data belongs to a subcritical space $H^s\times H^{s-1}$ for some $s>1$, \footnote{See however Tataru \cite{TataruRough} where a notion of finite energy solution was introduced. }  we shall need the following lemma when we deal with some initial data which is $C^{\infty}(R^2\backslash\{0\})$ but may fail to be in $H^s\times H^{s-1}$ globally for any $s>1$. 
\begin{lemma}\label{lm:exteriornice}
Suppose that $(u_0,\,u_1)\in C^{\infty}(R^2\backslash\{0\})$, and that $(u_0-u_{\infty},\,u_1)$ is compactly supported. Assume that 
\begin{equation}
\|(u_0,\,u_1)\|_{\HL}\leq \epsilon_{\ast}.
\end{equation}
Then there exists a unique smooth $u\in C^{\infty}((x,t):\,|x|>|t|\})$ such that $u$ solves the wave map equation in $\{(x,t):\,|x|>|t|\}$. Moreover
\begin{equation}\label{eq:exteriornice}
\lim_{t\to 0}\|\OR{u}(\cdot,\,t)-(u_0,\,u_1)\|_{\HL(|x|>|t|)}= 0.
\end{equation}
Similar results hold if we assume instead that $(u_0,\,u_1)\in H^{s}_{{\rm loc}}\times H^{s-1}_{{\rm loc}}(R^2\backslash\{0\})$, and in this case, $\OR{u}\in H^{s}_{{\rm loc}}\times H^{s-1}_{{\rm loc}}(|x|>|t|)$.
\end{lemma}

\smallskip
\noindent
{\it Proof.} We shall prove only the first part of the lemma. The proof of the second part is clear from the same argument. Let us firstly prove the existence of $u$ claimed in the lemma. For any $r>0$, since
$$\int_{B_r\backslash B_{\frac{r}{2}}}|\nabla u_0|^2+|u_1|^2\,dx\leq\epsilon_{\ast}^2,$$
we can find $\br\in \left(\frac{r}{2},\,r\right)$ with 
$$\int_{|x|=\br}|\spartial u_0|^2\,d\sigma\lesssim \frac{\epsilon^2_{\ast}}{\br}.$$
Denote
$$\overline{u}_0=\frac{1}{2\pi \br}\int_{\partial B_{\br}}u_0.$$
Then by Sobolev inequality, we get that
$$\left\|u_0-\overline{u}_0\right\|_{L^{\infty}(\partial B_{\br})}\lesssim \epsilon_{\ast}.$$
Thus from the fact that $|u_0|\equiv 1$, we see that $\left|\overline{u}_0\right|\gtrsim 1$. Take smooth cutoff function $\eta$ such that $\eta\equiv 1$ for $|x|\ge\br$ and $\eta\equiv 0$ for $|x|<\frac{\br}{2}$ with $|\nabla \eta|\lesssim (\br)^{-1}$.
Define
\begin{equation*}
\left(\widetilde{u}_0,\,\widetilde{u}_1\right)=\left\{\begin{array}{lr}
                                                       (u_0,\,u_1) & {\rm in}\,\,B_{\br}^c;\\
                           \left(P\left[\eta(r)(u_0(\br \theta)-\overline{u}_0)+\overline{u}_0\right],\,0\right) & {\rm in}\,\, B_{\br}.
                                                                            \end{array}\right.
\end{equation*}
Then
$$\left(\tu_0-u_{\infty},\,\tu_1\right)\in H^s\times H^{s-1},$$
for $s<\frac{3}{2}$, and direct computation shows that
$$\left\|\left(\widetilde{u}_0,\,\widetilde{u}_1\right)\right\|_{\HL}\lesssim \epsilon_{\ast}.$$
Note also that $\left(\widetilde{u}_0,\,\widetilde{u}_1\right)$ is smooth for $|x|>\br$.
 Hence by small data theory and Lemma \ref{lm:localsmoothness} the solution $\tu$ to the wave map equation with initial data $\left(\widetilde{u}_0,\,\widetilde{u}_1\right)$ is global, and is smooth in $|x|>\br+|t|$. By taking $r\to 0+$ and the finite speed of propagation, we see that 
$$u=\lim_{r\to 0+}\tu$$
exists in $|x|>|t|$ and is smooth. We now turn to the proof of (\ref{eq:exteriornice}). Let $\tu$ be the solution as before, corresponding to $\br$, then 
\begin{equation}\label{eq:equal90}
\tu\equiv u,\,\,\,{\rm for }\,\,\,|x|>\br+|t|,
\end{equation}
 and $\tu$ is continuous in $\HL$ for $t\in(0,\,1]$.  For any $\epsilon>0$, we can choose $\br$ sufficiently small, such that 
$$\|(u_0,\,u_1)\|_{\HL(B_{4\br})}<\epsilon.$$
Then by energy flux identity (say for $t>0$ and any $\epsilon>0$),
\begin{eqnarray*}
&&\int_{t+\epsilon<|x|<4\br-t}\nablaxtu(x,t)\,dx\\
&&\hspace{.3in}+\,\frac{1}{\sqrt{2}}\int_0^t\int_{|x|=4\br-t}\left(\frac{|\nabla u|^2}{2}+\frac{|\partial_tu|^2}{2}-\frac{x}{|x|}\cdot\nabla u\,\partial_tu\right)\,d\sigma ds\\
&&\hspace{.3in}+\,\frac{1}{\sqrt{2}}\int_0^t\int_{|x|=t+\epsilon}\left(\frac{|\nabla u|^2}{2}+\frac{|\partial_tu|^2}{2}+\frac{x}{|x|}\cdot\nabla u\,\partial_tu\right)\,d\sigma ds\\
&&\,\,=\,\int_{B_{4\br}\backslash B_{\epsilon}}\nablaxtu(x,0)\,dx,
\end{eqnarray*}
 we see that
\begin{equation}\label{eq:localsmall90}
\|\OR{u}(t)\|_{\HL\left(B_{2\br}\backslash B_{|t|}\right)}\leq \|(u_0,\,u_1)\|_{\HL\left(B_{4\br}\right)}<\epsilon,\,\,\,{\rm for}\,\,|t|<\br.
\end{equation}
Since $\tu$ is continuous in the energy space and $\OR{\tu}(x,0)=(u_0,\,u_1)(x)$ for $|x|>\br$, we see that for sufficiently small $t_{1}\in(0,\,\br)$ and $|t|<t_1$,
\begin{equation}\label{eq:tucontinuous}
\|\OR{\tu}(t)-(u_0,\,u_1)\|_{\HL(|x|>\br)}<\epsilon.
\end{equation}
Combining (\ref{eq:localsmall90}), (\ref{eq:tucontinuous}) and (\ref{eq:equal90}), we conclude that for $|t|\leq t_{1}$
\begin{eqnarray*}
&&\|\OR{u}(t)-(u_0,\,u_1)\|_{\HL(|x|>|t|)}\\
&&\leq \|\OR{u}(t)-(u_0,\,u_1)\|_{\HL(|x|>2\br)}+\|\OR{u}(t)-(u_0,\,u_1)\|_{\HL(|t|<|x|<2\br)}\\
&&\leq \|\OR{\tu}(t)-(u_0,\,u_1)\|_{\HL(|x|>2\br)}+2\epsilon\\
&&\leq 3\epsilon.
\end{eqnarray*}
Since $\epsilon>0$ is arbitrary, the lemma is proved.\\

 By finite speed of propagation and small data global existence, understanding the energy concentration is important for studying the dynamics of the wave maps. To measure the energy concentration, let us define for a wave map $u$ the ``{\it energy concentration radius}"
\begin{multline}\label{eq:concentrationradiusast}
\hspace{2in}r(\epsilon_{\ast},\,t):=\\
\inf\left\{r>0:\,{\rm there\,\,exists\,\,}\bx\,\,{\rm such\,\,that\,\,}\int_{B_r(\bx)}\nablaxtu(x,t)\,dx>\epsilon_{\ast}\right\}.
\end{multline}
We adopt the convention that if the set is empty, then the infimum is infinity.
The small energy global existence result, Theorem \ref{th:globalregularity}, and the finite speed of propagation imply that if wave map $u$ blows up at a finite time $T_+$, then $r(\epsilon_{\ast},t)\to 0+$ as $t\to T_+$. This is a very important piece of information that allows us to zoom in a small region near the blow up point and study the details of the blow up there. Unfortunately, knowing only that the energy concentrates in the small scales does not in itself allow one to ``extract" a nontrivial blow up profile in the limit, as we zooms in more and more. This is because a priori the energy can be concentrated in quite an arbitrary way, given that we do not (and it is probably not possible) to obtain control any higher order regularity beyond the energy when the time is close to the blow up time. To obtain a nontrivial blow up profile, the following result due to Sterbenz-Tataru \cite{TataruSterbenz} plays an essential role. \footnote{More precisely, this result is used to rule out the situation that 
all energy near the blow up point concentrates near the boundary of lightcone. The control inside the lightcone turns out to be quite favorable.}
\begin{theorem}\label{th:bubble}
There exists a function $\epsilon(E)$ with $0<\epsilon(E)\ll 1$ of the energy $E$ such that if $u$ is a classical solution to (\ref{Mainwavemap}) in $I\times R^2=[a,\,b]\times R^2$, with energy $E$ and 
\begin{equation}\label{eq:smalldispersion}
\sup_{t\in I}\,\sup_k\left\|(P_ku,\,2^{-k}P_k\partial_tu)(t)\right\|_{L^{\infty}\times L^{\infty}}<\epsilon(E),
\end{equation} 
then the energy concentration radius $r(\epsilon_{\ast},t)$ has a uniform lower bound on $I$:
\begin{equation}\label{eq:energyconcentrationradiuslowerbound}
\inf_{t\in I}\,r(\epsilon_{\ast},t)\ge r_0>0.
\end{equation}

\end{theorem}
\end{section}

\begin{section}{Channel of energy inequality for wave maps with small energy}
In this section, we prove the channel of energy inequality for small wave maps.
Let us begin with the following linear channel of energy inequality for outgoing waves, which is a slightly more quantitative two dimensional version of the channel of energy inequality that played a decisive role in \cite{DJKM}.
\begin{lemma}\label{lm:linearchannel}
Fix $\gamma\in(0,1)$. There exists $\mu=\mu(\gamma)>0$ sufficiently small such that the following statement is true. Let $v$ be a finite energy solution to the linear wave equation
$$\partial_{tt}v-\Delta v=0,\,\,{\rm in}\,\,R^2\times[0,\infty),$$
with initial data $(v_0,\,v_1)\in\HL$ satisfying
\begin{equation}\label{eq:outgoinglinearwave}
\|(v_0,\,v_1)\|_{\HL(B^c_{1+\mu}\cup B_{1-\mu})}+\|\spartial v_0\|_{L^2}+\|\partial_rv_0+v_1\|_{L^2}\leq \mu\|(v_0,\,v_1)\|_{\HL}.
\end{equation}
We also assume that $v_0\equiv v_{\infty}$ for some constant $v_{\infty}$ for large $x$.
Then for all $t\ge 0$, we have
\begin{equation}\label{eq:linearchannel}
\int_{|x|\ge \gamma+t}|\nabla_{x,t}v|^2(x,t)\,dx\ge \gamma\|(v_0,\,v_1)\|_{\HL}^2.
\end{equation}
\end{lemma}

\smallskip
\noindent
{\it Proof.} We can normalize the initial data so that $\|(v_0,\,v_1)\|_{\HL}=1$. 
Let $\alpha=\int_{\{1+\mu\leq |x|\leq 2(1+\mu)\}}v_0(x)\,dx$. By Poincar\'e inequality 
\begin{equation}
\label{eq:Poincare}
\int_{1+\mu\leq |x|\leq 2(1+\mu)} |v_0(x)-\alpha|^2\,dx\lesssim \int_{1+\mu\leq |x|\leq 2(1+\mu)} |\nabla v_0(x)|^2\,dx, 
\end{equation} 
where the implicit constant is independent of $\mu\leq 1$.

Take a non-negative radial $\eta\in C_c^{\infty}(R^2)$ with $\eta\equiv 1$ on $\overline{B_{1+\mu}}$ and ${\rm supp}\,\eta\Subset B_{1+\mu^{1/2}}$ satisfying $|\nabla \eta|\lesssim \mu^{-1/2}$. 
Define 
$$(\wt{v}_0,\wt{v}_1)=\eta(x)\left(v_0(x)-\alpha,v_1(x)\right).$$
Using \eqref{eq:Poincare}, the bound $|\nabla \eta|\lesssim \mu^{-1/2}$ and \eqref{eq:outgoinglinearwave}, we obtain:
\begin{equation}
\label{eq:veryclosechannel}
\left\|\left(\nabla(\wt{v}-v_0),\wt{v}_1-v_1)\right)\right\|_{L^2\times L^2}^2\lesssim \frac{1}{\mu}\int_{|x|\geq 1+\mu} |\nabla v_0|^2+\int_{|x|\geq 1+\mu} |u_1|^2\lesssim \mu.
\end{equation} 
By Sobolev and H\"older inequalities
\begin{gather*}
\||D|^{\frac 12} \wt{v}_0\|_{L^2}\lesssim \|\nabla \wt{v_0}\|_{L^{\frac 43}} \lesssim \|\nabla \wt{v}_0\|_{L^2(\{|x|\leq 1-\mu\})}+\mu^{\frac 18}\|\nabla \wt{v}_0\|_{L^2(\{1-\mu\leq |x|\leq 1+\mu^{1/2}\})}\lesssim \mu^{\frac 18}\\
\||D|^{-\frac 12} \wt{v}_1\|_{L^2}\lesssim \|v_1\|_{L^{\frac 43}}\lesssim \|v_1\|_{L^2(\{|x|\leq 1-\mu\})}+\mu^{\frac 18}\|v_1\|_{L^2(\{1-\mu\leq |x|\leq 1+\mu^{1/2}\})}\lesssim \mu^{\frac 18}.
\end{gather*}
By conservation of the $\dot{H}^{1/2}\times \dot{H}^{-1/2}$ norm for the linear wave equation, we obtain the for all $t\in \R$,
\begin{equation}
\label{eq:controlofl2}
\left|\int \wt{v}_t(x,t)\wt{v}(x,t)\,dx\right|\lesssim \big\| |D|^{1/2} \wt{v}\big\|_{L^2}\big\| |D|^{-1/2} \wt{v}_t\big\|_{L^2}\lesssim \mu^{1/4}.
\end{equation}

Let $\wt{v}$ be the solution to the linear wave equation with initial data $(\wt{v}_0,\,\wt{v}_1)$. 
By direct computation, we see that
\begin{equation}
\frac{d}{dt}\int_{R^2}-\wt{v}_t\left(x\cdot\nabla\wt{v}+\frac{1}{2}\wt{v}\right)(x,t)\,dx=\E(\wt{v}):=\E_0.
\end{equation}
Hence, by (\ref{eq:controlofl2}) and the outgoing condition (\ref{eq:outgoinglinearwave}), we get that
\begin{eqnarray*}
\int_{R^2}-\wt{v}_t\,x\cdot\nabla \wt{v}(x,t)\,dx&=&\E_0\, t+\int_{R^2}-\wt{v}_1\left(x\cdot\nabla\wt{v}_0+\frac{1}{2}\wt{v}_0\right)(x)\,dx\\
&&\hspace{.3in}+O(\mu^{1/4})\\
&=&\E_0\, (t+1)+O(\mu^{1/4}).
\end{eqnarray*}
On the other hand, by the finite speed of propagation, ${\rm supp}\,\wt{v}(\cdot,t)\Subset B_{1+\mu^{1/2}+t}$ for all $t\ge 0$, and thus
\begin{eqnarray*}
&&\int_{R^2}-\wt{v}_t\,x\cdot\nabla \wt{v}(x,t)\,dx\leq \int_{|x|>\gamma +t}(1+\mu^{1/2}+t)\left(\frac{\left|\wt{v}_t\right|^2}{2}+\frac{\left|\nabla \wt{v}\right|^2}{2}\right)(x,t)\,dx\\
&&\hspace{.3in} +\,(\gamma+t)\int_{|x|<\gamma+t}\left(\frac{\left|\wt{v}_t\right|^2}{2}+\frac{\left|\nabla \wt{v}\right|^2}{2}\right)(x,t)\,dx\\
&&=(\gamma+t)\E_0-(\gamma+t)\int_{|x|>\gamma+t}\left(\frac{\left|\wt{v}_t\right|^2}{2}+\frac{\left|\nabla \wt{v}\right|^2}{2}\right)(x,t)\,dx+\\
&&\hspace{.3in}+\,(1+\mu^{1/2}+t)\int_{|x|>\gamma+t}\left(\frac{\left|\wt{v}_t\right|^2}{2}+\frac{\left|\nabla \wt{v}\right|^2}{2}\right)(x,t)\,dx.
\end{eqnarray*}
Combining this and the above, we see that
\begin{eqnarray*}
&&(1+\mu^{1/2}-\gamma)\int_{|x|>\gamma+t}\left(\frac{\left|\wt{v}_t\right|^2}{2}+\frac{\left|\nabla \wt{v}\right|^2}{2}\right)(x,t)\,dx\\
&&\ge (1-\gamma)\E_0+O(\mu^{1/4}).
\end{eqnarray*}
By choosing $\mu$ sufficiently small, we obtain the channel of energy inequality for $\wt{v}$, and consequently also for $v$, by (\ref{eq:veryclosechannel}). \\

As mentioned in the introduction, one of the main goals of this paper is to extend the channel of energy arguments to the wave map setting. As a first step towards understanding the implications of the channel of energy property of linear wave equations on the wave maps, we prove the following result for small energy wave maps. The extension to large energy case seems to require nontrivial improvement in the perturbative techniques for the wave maps.
\begin{theorem}\label{th:channelofenergywavemap}
Fix $\beta\in (0,1)$. There exist a small  $\delta=\delta(\beta)>0$ and sufficiently small  $\epsilon_0=\epsilon_0(\beta)>0$,  such that if $u$ is a classical wave map with energy $\mathcal{E}(\OR{u})<\epsilon_0^2$ satisfying
\begin{equation}\label{eq:outgoingcondition}
\|(u_0,\,u_1)\|_{\HL\left(B^c_{1+\delta}\cup B_{1-\delta}\right)}+\|\spartial u_0\|_{L^2}+\|\partial_ru_0+u_1\|_{L^2}\leq \delta \|(u_0,\,u_1)\|_{\HL},
\end{equation}
then for all $t\ge 0$, we have
\begin{equation}\label{eq:smallchannel}
\int_{|x|>\beta+t}|\nabla_{x,t}u|^2(x,t)\,dx\ge \beta\,\|(u_0,\,u_1)\|_{\HL}^2.
\end{equation} 
\end{theorem}

\smallskip
\noindent
{\it Proof.} Denote $\epsilon:=\|(u_0,\,u_1)\|_{\HL}\lesssim \epsilon_0$. To apply Theorem \ref{th:globalregularity}, let us define the following frequency envelop
\begin{equation}
c_k:=\sup_{j\in\mathbb{Z}}2^{-\vartheta|k-j|}\|(P_ju_0,\,P_ju_1)\|_{\HL}.
\end{equation}
Then one can verify that $c=(c_k)$ is a frequency envelop and that $\left(\|P_k(u_0,u_1)\|_{\HL}\right)$ lies below it. In addition,
$$\|(c_k)\|_{\ell^2}\lesssim  \epsilon.$$
By Theorem \ref{th:globalregularity}, if $\epsilon_0$ is chosen sufficiently small, then the wave map $u$ is globally defined, and satisfies (\ref{eq:GlobalControl}). \\
Since the proof is a bit lengthy, we divide the arguments in several steps. \\

\smallskip
\noindent
{\bf Step 1: Reduction to proving channel of energy inequality for frequency pieces.}\\
In this step, our main goal is to show that there exists a set $\mathcal{K}$ of {\it good frequencies}, such that
\begin{equation}\label{eq:goodfrequencydominates}
\sum_{m\in\mathcal{K}}\|P_m(u_0,\,u_1)\|^2_{\HL}\ge (1-C\delta^{\frac{1}{12}})\initial^2,
\end{equation}
and that for any $m\in \mathcal{K}$, $2^m$ is {\it ``high frequency"}, and that it suffices to prove the channel of energy inequality for each $m\in \mathcal{K}$.\\

\noindent
{\it Substep (1): Control of the low frequency component.}\\
Fix $k_0$ large, whose precise value is to be determined below. We shall show that the total energy with frequency $\leq 2^{k_0}$ is small in a suitable sense. Assume firstly that $2^{-k_0}>C\delta$. Let us bound the low frequency energy of $(u_0,\,u_1)$, that is 
$$\left\|P_{\leq k_0}(u_0,\,u_1)\right\|_{\HL}.$$
 We can write
\begin{eqnarray*}
&&\nabla u_0=(\nabla u_0)\chi_{B_{1+\delta}^c\cup B_{1-\delta}}+(\nabla u_0)\chi_{B_{1+\delta}\backslash B_{1-\delta}},\\
&&u_1=u_1\chi_{B_{1+\delta}^c\cup B_{1-\delta}}+ u_1\chi_{B_{1+\delta}\backslash B_{1-\delta}}.
\end{eqnarray*}
By the assumption on $(u_0,\,u_1)$, 
\begin{eqnarray*}
\|(\nabla u_0)\chi_{B_{1+\delta}^c\cup B_{1-\delta}}\|_{L^2}&+&\|u_1\chi_{B_{1+\delta}^c\cup B_{1-\delta}}\|_{L^2}\\
&\lesssim&\delta \|(u_0,\,u_1)\|_{\HL}.
\end{eqnarray*}
Thus, 
\begin{eqnarray}
&&\left\|P_{<k_0}\left(\nabla u_0\chi_{B_{1+\delta}^c\cup B_{1-\delta}}\right)\right\|_{L^2}\nonumber\\
&&\quad\quad+\,\left\|P_{<k_0}\left(u_1\chi_{B_{1+\delta}^c\cup B_{1-\delta}}\right)\right\|_{L^2}\lesssim\delta \|(u_0,\,u_1)\|_{\HL}. \label{eq:low1}
\end{eqnarray}
Denote $f=(\nabla u_0)\chi_{B_{1+\delta}\backslash B_{1-\delta}}$. Then $f$ is compactly supported in $\overline{B_{1+\delta}}\backslash B_{1-\delta}$, and $$\|f\|_{L^2}\leq \|(u_0,\,u_1)\|_{\HL}.$$
By Bernstein's inequality, then Cauchy-Schwarz
$$\|P_{\leq k_0}f\|_{L^2}\lesssim 2^{k_0}\|f\|_{L^1}\lesssim 2^{k_0}\delta^{\frac 12}\|f\|_{L^2}\lesssim 2^{k_0}\delta^{\frac 12}\initial.$$

Choosing $2^{-k_0}\sim \delta^{\frac{1}{6}}$, then $\|P_{\leq k_0}f\|_{L^2}\lesssim \delta^{\frac{1}{3}}\initial$, that is,
\begin{equation}\label{eq:lowfrequencypotential}
\left\|P_{\leq k_0}\left[(\nabla u_0)\chi_{B_{1+\delta}\backslash B_{1-\delta}}\right]\right\|_{L^2}\lesssim \delta^{\frac{1}{3}}\initial.
\end{equation}
We can prove similarly that 
\begin{equation}\label{eq:lowfrequencytime}
\left\|P_{\leq k_0}\left[u_1\chi_{B_{1+\delta}\backslash B_{1-\delta}}\right]\right\|_{L^2}\lesssim \delta^{\frac{1}{3}}\initial.
\end{equation}
Combining (\ref{eq:low1}), (\ref{eq:lowfrequencypotential}) and (\ref{eq:lowfrequencytime}), we conclude that
\begin{equation}\label{eq:lowfreq}
\left\|P_{\leq k_0}(u_0,\,u_1)\right\|_{\HL}\lesssim \delta^{\frac{1}{3}}\initial.
\end{equation}
Thus the low frequency energy is small.\\

\medskip
\noindent
{\it Substep (2): persistence of condition (\ref{eq:outgoingcondition}) for most high frequencies.}\\
Let us now consider $P_k(\nabla u_0,\,u_1)$ for high frequency $2^{k}\ge2^{k_0}$. Fix small $\lambda>10\delta^{\frac{1}{6}}\sim 2^{-k_0}$, whose value is to be determined below.
Let us firstly bound $$\|P_k(\nabla u_0,\,u_1)\|_{L^2\left(B_{1+\lambda}^c\cup B_{1-\lambda}\right)}.$$
We can decompose as before
\begin{eqnarray*}
&&\nabla u_0=(\nabla u_0)\chi_{B_{1+\delta}^c\cup B_{1-\delta}}+(\nabla u_0)\chi_{B_{1+\delta}\backslash B_{1-\delta}},\\
&&u_1=u_1\chi_{B_{1+\delta}^c\cup B_{1-\delta}}+ u_1\chi_{B_{1+\delta}\backslash B_{1-\delta}}.
\end{eqnarray*}
Denote
$$\sigma_k:=\left\|P_k\left[(\nabla u_0,\,u_1)\chi_{B_{1+\delta}^c\cup B_{1-\delta}}\right]\right\|_{L^2\times L^2},$$
then it follows from (\ref{eq:outgoingcondition}) that
\begin{equation}
\sum_k\sigma_k^2\lesssim \delta^2\|(\nabla u_0,\,u_1)\|^2_{L^2\times L^2}.
\end{equation}
Now let us consider $P_k\left[(\nabla u_0,\,u_1)\chi_{B_{1+\delta}\backslash B_{1-\delta}}\right]$ for $x$ with $||x|-1|>\lambda$. Denote 
$$f:=(\nabla u_0) \chi_{B_{1+\delta}\backslash B_{1-\delta}},$$
then 
$$P_kf(x)=4^k\int_{R^2}\,\check{\Psi}(2^k(x-y))f(y)\,dy.$$
Since $f$ is supported in $1-\delta\leq |y|\leq 1+\delta$, and $||x|-1|>\lambda\gg \delta$, we get that
$$|P_kf(x)|\lesssim\frac{4^k}{(2^k||x|-1|)^M}\delta^{\frac{1}{2}}\initial.$$
Hence

\begin{gather*}
\|P_kf\|_{L^2\left(B_{1+\lambda}^c\cup B_{1-\lambda}\right)}^2\lesssim 4^{(2-M)k}\delta \initial^2\int_{\big||x|-1\big|\geq \lambda} \frac{1}{(|x|-1)^{2M}}\,dx\\
\lesssim 4^{(2-M)k}\delta\lambda^{-2M}\initial^2
\end{gather*}
Fix $M=3$. Then we conclude
\begin{equation}
\|P_kf\|_{L^2\left(B_{1+\lambda}^c\cup B_{1-\lambda}\right)}\leq 2^{-k}\delta^{\frac{1}{2}}\lambda^{-3}\initial.
\end{equation}
Take $\lambda=\delta^{\frac{1}{12}}$, then 
$$\|P_kf\|_{L^2\left(B_{1+\lambda}^c\cup B_{1-\lambda}\right)}\lesssim 2^{-k}\delta^{\frac{1}{4}}\initial,$$
that is, 
\begin{equation}\label{eq:controlpkf1}
\left\|P_k\left[(\nabla u_0)\chi_{B_{1+\delta}\backslash B_{1-\delta}}\right]\right\|_{L^2\left(B_{1+\lambda}^c\cup B_{1-\lambda}\right)}\lesssim 2^{-k}\delta^{\frac{1}{4}}\initial.
\end{equation}
Similarly, we can prove that
\begin{equation}\label{eq:controlpkf2}
\left\|P_k\left[u_1\chi_{B_{1+\delta}\backslash B_{1-\delta}}\right]\right\|_{L^2\left(B_{1+\lambda}^c\cup B_{1-\lambda}\right)}\lesssim 2^{-k}\delta^{\frac{1}{4}}\initial.
\end{equation}
Now let us control
$$\|\partial_rP_ku_0+P_ku_1\|_{L^2(B_{1+\lambda}\backslash B_{1-\lambda})}+\|\spartial P_ku_0\|_{L^2(B_{1+\lambda}\backslash B_{1-\lambda})}$$
for $k\ge k_0$ with $2^{-k_0}\sim \delta^{\frac{1}{6}}$. We have
\begin{eqnarray*}
&&\partial_r\int_{R^2}\,4^k\check{\Psi}(2^ky)u_0(x-y)\,dy=\int_{R^2}\,4^k\check{\Psi}(2^ky)\frac{x}{|x|}\cdot\nabla u_0(x-y)\,dy\\
&&=\int_{R^2}\,4^k\check{\Psi}(2^ky)\frac{x-y}{|x-y|}\cdot\nabla u_0(x-y)\,dy\\
&&\hspace{.4in}+\,\int_{R^2}\,4^k\check{\Psi}(2^ky)\left[\frac{x}{|x|}-\frac{x-y}{|x-y|}\right]\cdot\nabla u_0(x-y)\,dy\\
&&=I_k+II_k.
\end{eqnarray*}
Note that
$$I_k+P_ku_1=\int_{R^2}\,4^k\check{\Psi}(2^ky)(\partial_ru_0+u_1)(x-y)\,dy.$$
Thus
\begin{equation}\label{eq:control111}
\|I_k+P_ku_1\|_{L^2}\lesssim \|P_k(\partial_ru_0+u_1)\|_{L^2}.
\end{equation}
Note also that, for $x\in B_{1+\lambda}\backslash B_{1-\lambda}$,
$$\left|\nabla \frac{x}{|x|}\right|\lesssim 1,$$
thus,
\begin{eqnarray*}
|II_k|&\leq& \int_{R^2}\,4^k|\check{\Psi}|(2^ky)\left|\frac{x}{|x|}-\frac{x-y}{|x-y|}\right|\cdot|\nabla u_0(x-y)|\,dy\\
&\leq&\int_{|y|<2^{-\frac{k}{2}}}\,+\,\int_{|y|>2^{-\frac{k}{2}}}\\
&\lesssim&\int_{|y|<2^{-\frac{k}{2}}}2^{-\frac{k}{2}}4^k|\check{\Psi}|(2^ky)\cdot|\nabla u_0(x-y)|\,dy\\
&&\hspace{.3in}+\int_{|y|>2^{-\frac{k}{2}}}4^k\left|2^ky\right|^{-M}|\nabla u_0(x-y)|\,dy.
\end{eqnarray*}
Then simple computation shows that
\begin{eqnarray}
&&\|II_k\|_{L^2(B_{1+\lambda}\backslash B_{1-\lambda})}\nonumber\\
&&\lesssim \int_{|y|<2^{-\frac{k}{2}}}2^{-\frac{k}{2}}\,4^k|\check{\Psi}|(2^ky)\cdot\|\nabla u_0\|_{L^2}\,dy+\int_{|y|>2^{-\frac{k}{2}}}4^k\left|2^ky\right|^{-M}\|\nabla u_0\|_{L^2}\,dy\nonumber\\
&&\lesssim 2^{-\frac{k}{2}}\initial.\label{eq:control112}
\end{eqnarray}
Thus combining (\ref{eq:controlpkf1}), (\ref{eq:controlpkf2}), (\ref{eq:control111}) and (\ref{eq:control112}), we get that
\begin{equation}\label{eq:appout}
\|\partial_rP_ku_0+P_ku_1\|_{L^2(B_{1+\lambda}\backslash B_{1-\lambda})}\lesssim 2^{-\frac{k}{2}}\initial+\|P_k(\partial_ru_0+u_1)\|_{L^2}.
\end{equation}
The bound 
\begin{equation}\label{eq:appradial}
\|\spartial P_ku_0\|_{L^2(B_{1+\lambda}\backslash B_{1-\lambda})}\lesssim 2^{-\frac{k}{2}}\initial+\|P_k\,\spartial u_0\|_{L^2}
\end{equation} 
follows similarly from the previous arguments. \\

\medskip
\noindent
{\it Substep (3): Summary of estimates from substep (1) and substep (2) and the definition of good frequencies.}\\
From (\ref{eq:lowfreq}),(\ref{eq:controlpkf1}),(\ref{eq:controlpkf2}),(\ref{eq:appout}) and (\ref{eq:appradial}), we have, for $\delta^{\frac 16}\sim 2^{-k_0}$
\begin{eqnarray}
&(1)&\,\|P_{<k_0}(\nabla u_0,\,u_1)\|_{L^2\times L^2}\lesssim \delta^{\frac{1}{3}}\initial;\label{eq:(1)}\\
&(2)&\,\sum_{k\ge k_0}\left\|P_k\left[(\nabla u_0,\,u_1)\chi_{B_{1+\delta}^c\cup B_{1-\delta}}\right]\right\|^2_{L^2\times L^2}\lesssim \delta^2\initial^2;\label{eq:(1.1)}\\
&(3)&\,\sum_{k\ge k_0}\,\left\|P_k\left[(\nabla u_0,\,u_1)\chi_{B_{1+\delta}\backslash B_{1-\delta}}\right]\right\|_{L^2\times L^2(B_{1+\lambda}^c\cup B_{1-\lambda})}\nonumber\\
&&\,\,\lesssim \delta^{\frac{1}{4}}\initial;\label{eq:(2)}\\
&&\nonumber\\
&(4)&\|\partial_rP_ku_0+P_ku_1\|_{L^2(B_{1+\lambda}\backslash B_{1-\lambda})}+\|\spartial P_ku_0\|_{L^2(B_{1+\lambda}\backslash B_{1-\lambda})}\nonumber\\
&&\,\,\lesssim 2^{-\frac{k}{2}}\initial+\|P_k(\partial_ru_0+u_1)\|_{L^2}+\|P_k\,\spartial u_0\|_{L^2}.\label{eq:(3)}
\end{eqnarray}

By (\ref{eq:(1)}), we can focus on the high frequencies $2^{k}\ge 2^{k_0}$. Indeed, we have
\begin{equation}\label{eq:highfrequency}
\left\|P_{k\ge k_0}(\nabla u_0,\,u_1)\right\|_{L^2\times L^2}\ge \left(1-C\delta^{\frac{1}{4}}\right)\initial.
\end{equation}

By (\ref{eq:(1.1)}) and (\ref{eq:(2)}), we see that
\begin{eqnarray*}
&&\sum_{k\ge k_0}\|P_k(\nabla u_0,\,u_1)\|^2_{L^2\times L^2\left(B^c_{1+\lambda}\cup B_{1-\lambda}\right)}\\
&&\hspace{.7in}+\,P_k\left[(\nabla u_0,\,u_1)\chi_{B_{1+\delta}\backslash B_{1-\delta}}\right]\bigg\|^2_{L^2\times L^2(B^c_{1+\lambda}\cup B_{1-\lambda})}\\
&&\lesssim \sum_{k\ge k_0}\left\|P_k\left[(\nabla u_0,\,u_1)\chi_{B_{1+\delta}^c\cup B_{1-\delta}}\right]\right\|_{L^2\times L^2}^2\\
&&\hspace{.3in}+\,\sum_{k\ge k_0}\left\|P_k\left[(\nabla u_0,\,u_1)\chi_{B_{1+\delta}\backslash B_{1-\delta}}\right]\right\|^2_{L^2\times L^2\left(B^c_{1+\lambda}\cup B_{1-\lambda}\right)}\\
&&\lesssim\delta^{\frac{1}{2}}\initial^2.
\end{eqnarray*}

Then by the above calculation and (\ref{eq:(3)}), we get that
\begin{eqnarray*}
&&\sum_{k\ge k_0}\left[\|\partial_rP_ku_0+P_ku_1\|_{L^2}^2+\|\spartial P_ku_0\|_{L^2}^2\right]\\
&&\lesssim \sum_{k\ge k_0}\|P_k(\nabla u_0,\,u_1)\|^2_{L^2\times L^2\left(B_{1+\lambda}^c\cup B_{1-\lambda}\right)}+\\
&&\hspace{.3in}+\,\sum_{k\ge k_0}\left[\|\partial_rP_ku_0+P_ku_1\|^2_{L^2\left(B_{1+\lambda}\backslash B_{1-\lambda}\right)}+\|\spartial P_ku_0\|^2_{L^2\left(B_{1+\lambda}\backslash B_{1-\lambda}\right)}\right]\\
&&\lesssim\delta^{\frac{1}{2}}\initial^2+\sum_{k\ge k_0}2^{-k}\initial^2+\\
&&\hspace{.3in}+\,\sum_{k\ge k_0}\left(\|P_k(\partial_ru_0+u_1)\|^2_{L^2}+\|P_k\spartial u_0\|^2_{L^2}\right)\\
&&\lesssim \delta^{\frac{1}{6}}\initial^2.
\end{eqnarray*}

Hence, if we define the set
\begin{eqnarray*}
&&\mathcal{K}:=\bigg\{k\ge k_0:\,\|(P_ku_0,\,P_ku_1)\|_{\HL\left(B_{1+\lambda}^c\cup B_{1-\lambda}\right)}+\|\partial_rP_ku_0+P_ku_1\|_{L^2}\\
&&\hspace{.6in}+\,\|\spartial P_ku_0\|_{L^2}\leq \delta^{\frac{1}{100}}\|P_k(u_0,\,u_1)\|_{\HL}\bigg\},
\end{eqnarray*}
we can estimate that
\begin{eqnarray*}
&&\sum_{k\ge k_0,\,k\not\in\mathcal{K}}\|(P_ku_0,\,P_ku_1)\|_{\HL}^2\\
&&\lesssim \delta^{-\frac{1}{50}}\sum_{k\ge k_0,\,k\not\in\mathcal{K}}\bigg[\|(P_ku_0,\,P_ku_1)\|^2_{\HL\left(B_{1+\lambda}^c\cup B_{1-\lambda}\right)}\\
&&\hspace{1.2in}+\,\|\partial_rP_ku_0+P_ku_1\|^2_{L^2}+\|\spartial P_ku_0\|^2_{L^2}\bigg]\\
&&\lesssim \delta^{-\frac{1}{50}}\delta^{\frac{1}{6}}\initial^2\lesssim \delta^{\frac{1}{12}}\initial^2.
\end{eqnarray*}
Hence the total energy at frequencies $\sim2^k$ with $k\ge k_0,\,k\not\in \mathcal{K}$ is negligible, and
we will focus on the high frequency pieces $P_k(u_0,\,u_1)$ with $2^k\ge 2^{k_0}$ and $k\in\mathcal{K}$ below. \\

\medskip
\noindent
{\it Substep (4): Reduction to channel of energy inequality for frequencies in $\mathcal{K}$.}\\
Fix $m\in\mathcal{K}$, then 
\begin{multline}\label{eq:outgoingfork}
\|(P_mu_0,\,P_mu_1)\|_{\HL\left(B_{1+\lambda}^c\cup B_{1-\lambda}\right)}+\|\partial_rP_mu_0+P_mu_1\|_{L^2}+\|\spartial P_mu_0\|_{L^2}\\
\leq \delta^{\frac{1}{100}}\|P_m(u_0,\,u_1)\|_{\HL}.
\end{multline}
We claim that if we can show for each $m\in \mathcal{K}$ that
\begin{equation}\label{eq:channelform}
\int_{|x|\ge \frac{1+\beta}{2}+t}\,|\nabla_{x,t}P_mu|^2(x,t)\,dx\ge \frac{1+\beta}{2}\|P_m(u_0,\,u_1)\|_{\HL}^2-C\epsilon^2c_m^2,
\end{equation}
for all $t\ge 0$, then we will be done. Indeed, write for each $t\ge 0$,
$$P_m\nabla_{x,t}u=P_m\left[(\nabla_{x,t}u)\chi_{|x|>\beta+t}\right]+P_m\left[(\nabla_{x,t}u)\chi_{|x|\leq\beta+t}\right].$$
We can estimate, for $|x|>\frac{\beta+1}{2}+t$, that
\begin{eqnarray*}
&&\left|P_m\left[(\nabla_{x,t}u)\chi_{|y|\leq\beta+t}\right](x)\right|\\
&&\leq 4^m\int_{|x-y|\leq \beta+t}\big|\check{\Psi}(2^my)\big|\big|\nabla_{x,t}u(x-y,t)\big|\,dy\\
&&\leq 4^m\int_{|y|>\frac{1-\beta}{2}}\big|\check{\Psi}(2^my)\big||\nabla_{x,t}u(x-y,t)|\,dy\\
&&\leq 4^m\int_{|y|>\frac{1-\beta}{2}}\big|2^my\big|^{-M}|\nabla_{x,t}u(x-y,t)|\,dy.
\end{eqnarray*}
Thus,
\begin{eqnarray*}
&&\left\|P_m\left[|\nabla_{x,t}u|\chi_{|x|\leq \beta+t}\right]\right\|_{L^2\left(|x|>\frac{\beta+1}{2}+t\right)}\\
&&\lesssim 4^m\int_{|y|>\frac{1-\beta}{2}}\left|2^my\right|^{-M}\,dy\,\initial\\
&&\lesssim C(\beta)\,2^{-(M-2)m}\initial.
\end{eqnarray*}
Consequently, we get that
\begin{eqnarray*}
&&\sum_{m\ge k_0}\|P_m\nabla_{x,t}u\|^2_{L^2\left(|x|>\frac{1+\beta}{2}+t\right)}\\
&&\leq \left(1+\delta^{\frac{1}{50}}\right)\sum_{m\ge k_0}\left\|P_m\left[(\nabla_{x,t}u)\chi_{|y|\ge \beta+t}\right]\right\|^2_{L^2}+\\
&&\hspace{.5in}+\,2\delta^{-\frac{1}{50}}\sum_{m\ge k_0}\left\|P_m\left[(\nabla_{x,t}u)\chi_{|y|\leq \beta+t}\right]\right\|_{L^2\left(|x|>\frac{\beta+1}{2}+t\right)}^2\\
&&\leq \left(1+\delta^{\frac{1}{50}}\right)\int_{|x|>\beta+t}|\nabla_{x,t}u|^2(x,t)\,dx+\\
&&\hspace{.6in}+\,C(\beta)\delta^{-\frac{1}{50}}\sum_{m\ge k_0}2^{-2(M-2)m}\initial^2\\
&&\leq \left(1+\delta^{\frac{1}{50}}\right)\int_{|x|>\beta+t}|\nabla_{x,t}u|^2(x,t)\,dx+C(\beta)\delta^{-\frac{1}{50}}4^{-(M-2)k_0}\initial^2\\
&&\leq \left(1+\delta^{\frac{1}{50}}\right)\int_{|x|>\beta+t}|\nabla_{x,t}u|^2(x,t)\,dx+C(\beta)\delta^{\frac{1}{3}}\initial^2,
\end{eqnarray*}
if we choose $M=4$.
Therefore if (\ref{eq:channelform}) holds, then by the choice of $\mathcal{K}$, (\ref{eq:highfrequency}), $\|(c_k)\|_{l^2}\lesssim \epsilon$, and the above calculation, we see that
\begin{eqnarray*}
&&(1-C\delta^{\frac{1}{12}}-C\epsilon^2)\frac{1+\beta}{2}\initial^2\\
&&\leq\sum_{m\in\mathcal{K}}\left(\frac{1+\beta}{2}\|P_m(u_0,\,u_1)\|^2_{\HL}-C^2\epsilon^2c_m^2\right)\\
&&\leq \sum_{m\in\mathcal{K}}\|P_m\nabla_{x,t}u\|^2_{L^2\left(|x|>\frac{1+\beta}{2}+t\right)}\\
&&\leq\left(1+\delta^{\frac{1}{50}}\right)\int_{|x|>\beta+t}|\nabla_{x,t}u|^2(x,t)\,dx+C(\beta)\delta^{\frac{1}{3}}\initial^2.
\end{eqnarray*}
The channel of energy inequality (\ref{eq:smallchannel}) follows if $\delta=\delta(\beta)$ and $\epsilon=\epsilon_0(\beta)$ are taken sufficiently small. Our goal is thus reduced to proving (\ref{eq:channelform}). \\

\medskip
\noindent
{\bf Step 2: Control of of the perturbative part of the nonlinearity.}\\
 It is proved in \cite{DJKM} that (\ref{eq:channelform}) holds for solution to the linear wave equation with this type of outgoing initial data for dimension $\ge 3$, although the results we need here are more quantitative, see Lemma \ref{lm:linearchannel} above.  Ideally one would like to say that the nonlinearity is negligible as we have small solutions. However, as is now well known, even in small energy case, the nonlinearity for the wave map equation can not be treated entirely perturbatively. Rather, we need to perform a gauge transform to modify the nonlinearity so that it becomes perturbative. Thus it is important to understand how the Gauge transform affects the channel of energy inequality. The arguments we use here are mostly from Tao \cite{TaoSmallEnergy} and Tataru \cite{Tataru3}. We shall present the details of the proof below, partly for the convenience of the reader, and partly as those works did not explicitly quantify the nonlinear effects (which are implicit in the proofs). In this 
step 
however, we shall firstly control the part of the nonlinearity that is perturbative.\\

Let $$\psi:=P_mu.$$ Then $\psi$ verifies
\begin{equation}\label{eq:forpsi}
\left\{\begin{array}{rll}
         \partial_{tt}\psi-\Delta\psi&=&P_m\left(u\,\partial^{\alpha}u^{\dagger}\partial_{\alpha}u\right)\\
           \OR{\psi}(0)&=&(P_mu_0,\,P_mu_1).
         \end{array}\right.
\end{equation}

Let us rewrite the nonlinearity $P_m\left(u\,\partial^{\alpha}u^{\dagger}\partial_{\alpha}u\right)$ as
\begin{eqnarray*} 
&&P_m\left(u\,\partial^{\alpha}u^{\dagger}\partial_{\alpha}u\right)\\
&&=P_m\left(u_{\ge m-10}\,\partial^{\alpha}u^{\dagger}\partial_{\alpha}u\right)\\
&&\hspace{.3in}+\,P_m\left(u_{<m-10}\,\partial^{\alpha}u^{\dagger}_{>m+10}\partial_{\alpha}u\right)\\
&&\hspace{.3in}+\,P_m\left(u_{<m-10}\,\partial^{\alpha}u^{\dagger}_{m-10\leq\cdot\leq m+10}\partial_{\alpha}u_{\ge m-10}\right)\\
&&\hspace{.3in}+\,P_m\left(u_{<m-10}\,\partial^{\alpha}u^{\dagger}_{m-10\leq \cdot\leq m+10}\partial_{\alpha}u_{< m-10}\right)\\
&&\hspace{.3in}+\,P_m\left(u_{<m-10}\,\partial^{\alpha}u^{\dagger}_{<m-10}\partial_{\alpha}u_{> m+10}\right)\\
&&\hspace{.3in}+\,P_m\left(u_{<m-10}\,\partial^{\alpha}u^{\dagger}_{<m-10}\partial_{\alpha}u_{m-10\leq\cdot\leq m+10}\right)\\
&&\hspace{.3in}+\,P_m\left(u_{<m-10}\,\partial^{\alpha}u^{\dagger}_{<m-10}\partial_{\alpha}u_{<m-10}\right)\\
&&=I_1+I_2+I_3+I_4+I_5+I_6+I_7.
\end{eqnarray*}
Denote 
$$\epsilon:=\initial\leq \epsilon_0.$$
We firstly peel off the perturbative part of the nonlinearity. We shall call $h$ {\it disposable} if
$$\sup_{m'=m+O(1)}\|P_{m'}h\|_{N[m']}\lesssim \epsilon\,c_m.$$
Here the $O(1)$ term is a number of size $\sim 10$. The main use of this term is to deal with some technical ``frequency leakage" issues. \footnote{On a technical level, to apply the estimates from Theorem \ref{th:propertySN}, we need the right hand sides to carry the frequency localization operator $P_k$. } We shall call $h$ {\it disposable in the generalized sense} if there exists a sequence of disposable $h_k$ with $h_k\to h$ in the sense of distributions. Note that the notion of being disposable and that of being disposable in the generalized sense  are not the same, due to the technical issue with the space $N[m]$, see e.g., page 324 of \cite{TataruRough} for more discussions.

Note that $I_4=I_6$. Furthermore, analyzing the support of the trilinear expressions in frequencies, we obtain that $I_5=I_7=0$.
We claim that $I_1,\,I_2,\,I_3$ are disposable, that is, 
\begin{claim}\label{claim:Ijareperturbative}
For $j=1,\,2,\,3$ we have
\begin{equation}\label{eq:perturbativenonlinearity}
\sup_{m'=m+O(1)}\left\|P_{m'}I_j\right\|_{N[m']}\lesssim \epsilon c_m.
\end{equation}
\end{claim}
We also claim that
\begin{claim}\label{claim:commutingperturbative}
For $m'=m+O(1)$,
\begin{multline}\label{eq:commutatvieerror}
\left\|P_{m'}\left[P_m\left(u_{<m-10}\,\partial^{\alpha}u^{\dagger}_{m-10<\cdot<m+10}\partial_{\alpha}u_{<m-10}\right)-u_{<m-10}\,\partial^{\alpha}\psi^{\dagger}\partial_{\alpha}u_{<m-10}\right]\right\|_{N[m']}\\
\lesssim \epsilon c_m,
\end{multline}
where $\psi$ is defined in (\ref{eq:forpsi}); similarly,
\begin{multline}\label{eq:commutatvieerror2}
\left\|P_{m'}\left[P_m\left(\partial_{\alpha}u_{<m-10}\,u^{\dagger}_{<m-10}\,\partial^{\alpha}u_{m-10<\cdot<m+10}\right)-\partial_{\alpha}u_{<m-10}\,u^{\dagger}_{<m-10}\,\partial^{\alpha}\psi\right]\right\|_{N[m']}\\
\lesssim \epsilon c_m.
\end{multline}
\end{claim}
We postpone the proof of Claim \ref{claim:Ijareperturbative} and Claim \ref{claim:commutingperturbative} to the end of this section.

Hence, by (\ref{eq:commutatvieerror}) we can rewrite the equation for $\psi$ as
\begin{equation}\label{eq:intermediateform}
\partial_{tt}\psi-\Delta\psi=\widetilde{f}+2u_{<m-10}\,\partial_{\alpha}u^{\dagger}_{<m-10}\partial^{\alpha}\psi,
\end{equation}
where
$$\sup_{m'=m+O(1)}\left\|P_{m'}\widetilde{f}\right\|_{N[m']}\lesssim \epsilon c_m.$$
Let us note the relation 
$$u^{\dagger}\,\partial_{\alpha}u=0.$$
It follows that
\begin{eqnarray*}
0&=&P_m\left(\partial^{\alpha}u_{<m-10}\ud\,\partial_{\alpha}u\right)\\
&=&P_m\left(\partial^{\alpha}u_{<m-10}\ud_{\ge m-10}\partial_{\alpha}u\right)+P_m\left(\partial^{\alpha}u_{<m-10}\ud_{< m-10}\partial_{\alpha}u_{m-10\leq\cdot\leq m+10}\right)\\
&=&I+II.
\end{eqnarray*}
We can estimate the $I$ term, by using the trilinear estimate, as for $m'=m+O(1)$
\begin{eqnarray*}
&&\left\|P_{m'}\left(\partial^{\alpha}u_{<m-10}\ud_{\ge m-10}\partial_{\alpha}u\right)\right\|_{N[m']}\\
&&\lesssim \sum_{k_1<m-10,\,m-10\leq k_2\leq m+10,\,k_3\leq m+O(1)}\left\|P_{m'}\left(\partial^{\alpha}u_{k_1}\ud_{k_2}\,\partial_{\alpha}u_{k_3}\right)\right\|_{N[m']}+\\
&&\hspace{.5in}+\,\sum_{k_1<m-10,\,k_2>m+10,\,k_3=k_2+O(1)}\left\|P_{m'}\left(\partial^{\alpha}u_{k_1}\ud_{k_2}\,\partial_{\alpha}u_{k_3}\right)\right\|_{N[m']}\\
&&\lesssim \sum_{k_1<m-10,\,k_3\leq m+O(1)} 2^{-\kappa (m-\min\{k_1,\,k_3\})}c_{k_1}c_mc_{k_3}+\\
&&\hspace{.5in}+\,\sum_{k_1<m-10,\,k_2>m+10}2^{-\kappa(k_2-m)}2^{-\kappa(k_2-k_1)}c_{k_1}c_{k_2}^2\\
&&\lesssim\epsilon^2 c_m
\end{eqnarray*}

Consequently, by the boundedness of $P_{m}$ in $S[m']$, we see that
\begin{equation*}
\sup_{m'=m+O(1)}\|P_{m'}II\|_{N[m']}\lesssim \epsilon^2 c_m,
\end{equation*}
and thus $II$ is disposable.
Thus by (\ref{eq:commutatvieerror2}) we can rewrite the equation for $\psi$ as
\begin{equation}\label{eq:mainequationforpsi}
\partial_{tt}\psi-\Delta\psi=f+2\left(u_{<m-10}\,\partial^{\alpha}u^{\dagger}_{<m-10}-\partial^{\alpha}u_{<m-10}\,u^{\dagger}_{<m-10}\right)\partial_{\alpha}\psi,
\end{equation}
where $f$ satisfies 
\begin{equation}
\sup_{m'=m+O(1)}\|f\|_{N[m']}\lesssim \epsilon c_m.
\end{equation}
\\

\medskip
\noindent
{\bf Step 3: Construction of the micro-local gauge.}\\

To deal with the non-perturbative part of the nonlinearity, we need to use the idea of Tao \cite{TaoSmallEnergy}.\\

We have
\begin{equation}\label{eq:modified}
\partial_{tt}\psi-\Delta\psi=f+2\left(u_{<m-10}\,\partial^{\alpha}u^{\dagger}_{<m-10}-\partial^{\alpha}u_{<m-10}\,u^{\dagger}_{<m-10}\right)\partial_{\alpha}\psi,
\end{equation}
where $f$ satisfies
\begin{equation}\label{eq:disposablef}
\sup_{m'=m+O(1)}\|f\|_{N[m']}\lesssim \epsilon c_m.
\end{equation}
Let $w=U_{<m-10}\psi$ for some matrix $U_{<m-10}$ to be determined below, then (\ref{eq:modified}) implies that
\begin{eqnarray*}
&&-\partial^{\alpha}\partial_{\alpha}w=-\partial^{\alpha}\partial_{\alpha}U_{<m-10}\,\psi-2\partial^{\alpha}U_{<m-10}\partial_{\alpha}\psi\\
&&\hspace{.5in}+\,U_{<m-10}\left[f+2\left(u_{<m-10}\,\partial^{\alpha}u^{\dagger}_{<m-10}-\partial^{\alpha}u_{<m-10}\,u^{\dagger}_{<m-10}\right)\partial_{\alpha}\psi\right]\\
&&=\left(\Box U_{<m-10}\right)\,\psi+U_{<m-10}f+\\
&&\hspace{.3in}+\,2\left[U_{<m-10}\left(u_{<m-10}\,\partial^{\alpha}u^{\dagger}_{<m-10}-\partial^{\alpha}u_{<m-10}\,u^{\dagger}_{<m-10}\right)-\partial^{\alpha} U_{<m-10}\right]\partial_{\alpha}\psi.
\end{eqnarray*}
Then
\begin{eqnarray}
&&\partial_{tt}w-\Delta w=\left(\Box U_{<m-10}\right)\psi+U_{<m-10}f\label{eq:modifiedwithg}\\
&&+\,2\left[U_{<m-10}\left(u_{<m-10}\,\partial^{\alpha}u^{\dagger}_{<m-10}-\partial^{\alpha}u_{<m-10}\,u^{\dagger}_{<m-10}\right)-\partial^{\alpha} U_{<m-10}\right]\partial_{\alpha}\psi.\nonumber
\end{eqnarray}
Ideally we would like to choose $U_{<m-10}$ so that
$$\partial^{\alpha}U_{<m-10}=U_{<m-10}\left(u_{<m-10}\partial^{\alpha}u^{\dagger}_{<m-10}-\partial^{\alpha}u_{<m-10}\,u^{\dagger}_{<m-10}\right),$$
for all $\alpha$, then the term on the right hand side of (\ref{eq:modifiedwithg}) containing $\partial_{\alpha}\psi$ would be eliminated, and we would be in a truly semilinear case. However this is impossible due to compatibility issues, see the discussions in \cite{TaoSmallEnergy}. Instead we will follow Tataru's modification of Tao's idea in \cite{TataruRough} to construct a micro-local approximate solution. \\

Fix large $N>1$. Define inductively
\begin{eqnarray*}
&&U^N_{-N}=I;\\
&&U^N_k=U^N_{<k-10}\left(u_{<k-10}u^{\dagger}_k-u_ku^{\dagger}_{<k-10}\right),
\end{eqnarray*}
where $U^N_{<k-10}=\sum\limits_{-N<j<k-10}U^N_j+I$ if $k>-N+11$ and $U^N_{<k-10}=I$ otherwise. In the end we will pass $N\to\infty$, but we need to obtain uniform in $N$ estimates for $U^N_k$ in order to do that. We claim the following properties for $U^N_k$ and $U^N_{<k}$ with $-N<k\leq N$:
\begin{claim}\label{cl:uniformU}
For $\epsilon$ sufficiently small,
\begin{eqnarray}
&&U^N_k\,\,{\rm has\,\,frequency\,\,support\,\,} 2^{k-2}\leq |\xi|\leq 2^{k+2};\\
&&\sup_{k'=k+O(1)}\left\|P_{k'}U^N_k\right\|_{S[k']}\lesssim c_k;\label{eq:bound_ck}\\
&&\left\|U^N_{<k}\left(U^N_{<k}\right)^{\dagger}-I\right\|_{S(c)}\lesssim \sqrt{\epsilon}.\label{eq:orthSc}
\end{eqnarray}
\end{claim}
We shall prove the claim inductively. For $k=-N+1$, the claim follows from the property that by Theorem \ref{th:globalregularity}
$$\|u\|_{S(c)}\lesssim 1.$$ Suppose the claim is true up to $k-1$, let us prove it holds also for $k$. A crucial point is the following important algebraic identity:
\begin{eqnarray}
&&U^N_k\left(U^N_{<k-10}\right)^{\dagger}+U^N_{<k-10}(U^N_k)^{\dagger}\\
&&=U^N_{<k-10}\left(u_{<k-10}u^{\dagger}_k-u_ku^{\dagger}_{<k-10}\right)\left(U^N_{<k-10}\right)^{\dagger}+\\
&&\hspace{.3in}+\,U^N_{<k-10}\left(u_ku^{\dagger}_{<k-10}-u_{<k-10}u^{\dagger}_k\right)\left(U^N_{<k-10}\right)^{\dagger}\\
&&=0.
\end{eqnarray}
We also note that $U^N_j$ is anti-symmetric if $-N<j\leq-N+11$, which is an easy consequence of the definition of $U^N_j$.\\
Thus by the anti-symmetry of $U^N_j$ for $-N<j\leq -N+11$, we get that
\begin{eqnarray*}
&&U^N_{<k}\left(U^N_{<k}\right)^{\dagger}=\left(\sum_{-N\leq j<k}U^N_j\right)\left(\sum_{-N\leq j<k}\left(U^N_j\right)^{\dagger}\right)\\
&&=\sum_{-N\leq j<j'-10<j'<k}U^N_j\left(U^N_{j'}\right)^{\dagger}+\sum_{-N\leq j'<j-10<j<k}U^N_j\left(U^N_{j'}\right)^{\dagger}+\\
&&\\
&&\hspace{.3in}+\,\sum_{|j-j'|\leq 10, \,-N< j,\,j'<k} U^N_j\left(U^N_{j'}\right)^{\dagger}+\sum_{-N<j\leq -N+10}\left[\left(U^N_{j}\right)^{\dagger}+U^N_j\right]+I\\
&&\\
&&=\sum_{-N+10< j'<k}U^N_{<j'-10}\left(U^N_{j'}\right)^{\dagger}+\sum_{-N+10< j<k}U^N_{j}\left(U^N_{<j-10}\right)^{\dagger}+I+\\
&&\hspace{.5in}+\,\sum_{|j-j'|\leq 10, \,-N< j,\,j'<k} U^N_j\left(U^N_{j'}\right)^{\dagger}\\
&&=I+\sum_{-N<j,\,j'<k,\,|j-j'|\leq 10}U^N_j\left(U^N_{j'}\right)^{\dagger}.
\end{eqnarray*}
Simplifying the above, we get that
$$U^N_{<k}\left(U^N_{<k}\right)^{\dagger}-I=\sum_{-N<j,\,j'<k,\,|j-j'|\leq 10}U^N_j\left(U^N_{j'}\right)^{\dagger}.$$
Hence by (\ref{eq:bound_ck}) from induction,
\begin{eqnarray*}
&&\left\|U^N_{<k}\left(U^N_{<k}\right)^{\dagger}-I\right\|_{L^{\infty}}\\
&&\lesssim \sum_{-N<j,\,j'<k,\,|j-j'|\leq 10}\left\|U^N_j\right\|_{L^{\infty}}\left\|U^N_{j'}\right\|_{L^{\infty}}\\
&&\lesssim \sum_{-N<j,\,j'<k,\,|j-j'|\leq 10}\,\,\sum_{j_1=j+O(1),\,j_2=j'+O(1)}\left\|P_{j_1}U^N_j\right\|_{S[j_1]}\left\|P_{j_2}U^N_{j'}\right\|_{S[j_2]}\\
&&\lesssim \sum_{-N<j<k}c_j^2\lesssim \epsilon^2.
\end{eqnarray*}
In the second inequality above, we used the fact that $U^N_j=\sum_{j_1=j+O(1)}P_{j_1}U_j^N$
which follows from the frequency support property of $U^N_j$. We shall use this trick often, as a replacement of bound on $\|U_j^N\|_{S[j]}$ which we do not have. Below we will omit the routine details when we use the same trick. In particular, combining the above with the induction bound (\ref{eq:bound_ck}), we see that $\|U^N_{<k}\|_{S(1)}\leq C$ for some universal constant (by choosing $\epsilon$ small).\\
Similarly, for each $k'<k+O(1)$, by the property of $S[k]$ spaces and induction,
\begin{eqnarray*}
&&\left\|P_{k'}\left[U_{<k}^N\left(U^N_{<k}\right)^{\dagger}\right]\right\|_{S[k']}\\
&&\lesssim \sum_{-N<j,\,j'<k,\,|j-j'|\leq 10}\left\|P_{k'}\left[U^N_j\left(U^N_{j'}\right)^{\dagger}\right]\right\|_{S[k']}\\
&&\lesssim\sum_{O(1)+k'<j<k}2^{-\kappa (j-k')_{+}}c_j^2\\
&&\lesssim \sum_{O(1)+k'<j<k}2^{-(\kappa-\vartheta) (j-k')_{+}}c_jc_{k'}\lesssim \epsilon\,c_{k'}.
\end{eqnarray*}
Combining the above two estimates, (\ref{eq:orthSc}) follows.  \\
The estimate for $\sup_{k'=k+O(1)}\left\|P_{k'}U^N_k\right\|_{S[k']}$ then follows from the definition and the fact that $\|u_{<k-10}\|_{S(c)}$, $\|U^N_{<k-10}\|_{S(1)}$ are universally bounded. The support property is obvious. \\

Using these uniform estimates, we can pass $N\to\infty$, and obtain a limit along a subsequence of $N$, so that 
$$U_k:=\lim_{N_i\to\infty}U^{N_i}_{k}, \,\,\,\,U_{<k}:=\lim_{N_i\to\infty}U^{N_i}_{<k},$$
exist in the sense of distributions, for each $k$. Since $U^N_k$ are frequency localized and have bounded overlap in frequency support, we can conclude that
\begin{equation}\label{eq:iterativeU}
U_{<k}=\sum_{k'<k}U_{k'}+I,\,\,{\rm and}\,\,U_k=U_{<k-10}(u_{<k-10}\ud_k-u_k\ud_{<k-10}).
\end{equation}
In addition, $U_k$, $U_{<k}$ satisfies the same estimates claimed for $U^N_k,\,U^N_{<k}$ above. As a consequence, we have 
\begin{equation}\label{eq:boundforUk}
\sup_{k'=k+O(1)}\left\|P_{k'}U_k\right\|_{S[k']}\lesssim c_k,\,\,\,\,{\rm and}\,\,\left\|U_{<k}\right\|_{S(c)}\lesssim 1.
\end{equation}
This is a direct consequence of the property of $S[k]$ under weak convergence, see the remark below Theorem \ref{th:propertySN}. \\

\medskip
\noindent
{\bf Step 4: Control of the nonlinearity after applying the gauge transform.}\\
We shall show that the terms on the right hand size of (\ref{eq:modifiedwithg}) are all disposable.\\

\noindent
{\it Substep (1): the terms involving $\Box U_k$.}\\
To control the terms $\left(\Box U_{<m-10}\right)\psi$, we need to control $\Box U^N_{<m-10}$ uniformly for all large $N$. By definition, 
\begin{eqnarray*}
&&\Box U^N_k=\left(\Box U^N_{<k-10}\right)\left(u_{<k-10}u^{\dagger}_k-u_ku^{\dagger}_{<k-10}\right)\\
&&\hspace{.4in}-\,2\partial^{\alpha}U^N_{<k-10}\partial_{\alpha}\left(u_{<k-10}u^{\dagger}_k-u_ku^{\dagger}_{<k-10}\right)\\
&&\hspace{.4in}+\,U^N_{<k-10}\left(\Box u_{<k-10}u^{\dagger}_k+u_{<k-10}\Box u^{\dagger}_k-\Box u_k\,u^{\dagger}_{<k-10}-u_k\Box u_{<k-10}^{\dagger}\right)\\
&&\hspace{.4in}+\,2U^N_{<k-10}\left(\partial^{\alpha}u_k\partial_{\alpha} \ud_{<k-10}-\partial^{\alpha}u_{<k-10}\partial_{\alpha}u_k^{\dagger}\right)=I+II+III+IV.
\end{eqnarray*}
We claim that for $\nu=\frac{\kappa}{32}$, and uniformly for all large $N$
\begin{claim}\label{claim:boxUterm}
\begin{equation}\label{eq:boxUestimate}
\sup_{j'=j+O(1)}\left\|P_{j'}\left(\Box U^N_k\phi\right)\right\|_{N[j']}\lesssim 2^{-\nu(j-k)}c_k\|\phi\|_{S[j]},
\end{equation}
for all $\phi$ with frequency support in $2^{j-5/2}\leq \cdot\leq 2^{j+5/2}$ and $k<j-7$. 
\end{claim}
Assuming this claim for a moment, then we can estimate for $m'=m+O(1)$
\begin{eqnarray*}
&&\left\|P_{m'}\left[\Box U^N_{<m-10}\,\psi\right]\right\|_{N[m']}\\
&&\lesssim \sum_{k<m-10}\left\|P_{m'}\left[\Box U^N_k\,\psi\right]\right\|_{N[m']}\\
&&\lesssim \sum_{k<m-10}2^{-\nu(m-k)}c_kc_m\lesssim \epsilon\, c_m,
\end{eqnarray*}
and thus the first term on the right hand side of (\ref{eq:modifiedwithg}) is disposable in the generalized sense.\\

We shall prove (\ref{eq:boxUestimate}) inductively. It is clear that (\ref{eq:boxUestimate}) holds for $k=-N$. Suppose (\ref{eq:boxUestimate}) holds for $k'<k$, let us prove that it holds for $k$. The bound for $I$ term: 
\begin{eqnarray*}
&&\left\|P_{j'}\left[\Box U^N_{<k-10}\left(u_{<k-10}u^{\dagger}_k-u_k\,u^{\dagger}_{<k-10}\right)\phi\right]\right\|_{N[j']}\\
&&\lesssim\sum _{k'<k-10,\,|j-j''|\leq 3}\left\|P_{j'}\left[\Box U^N_{k'}\,P_{j''}\left\{\left(u_{<k-10}u^{\dagger}_k-u_k\,u^{\dagger}_{<k-10}\right)\phi\right\}\right]\right\|_{N[j']}\\
&&\lesssim \sum_{k'<k-10}2^{-\nu (j-k')}c_{k'}c_k\|\phi\|_{S[j]}\\
&&\lesssim 2^{-\nu(j-k)}\epsilon\,c_k\|\phi\|_{S[j]}
\end{eqnarray*}
follows from the inductive hypothesis and the property of $S[k]$ spaces. The projection $P_{j''}$ was used to deal with the frequency leakage, which is a minor technical issue.\\

Let us consider the $II$ term $\partial^{\alpha}U^N_{<k-10}\partial_{\alpha}\left(u_{<k-10}\ud_k-u_k\ud_{<k-10}\right)\phi$. By (\ref{eq:boundforUk}) and the trilinear estimate, we have
\begin{eqnarray*}
&&\left\|P_{j'}\left[\partial^{\alpha}U^N_{<k-10}\partial_{\alpha}\left(u_{<k-10}\ud_k-u_k\ud_{<k-10}\right)\phi\right]\right\|_{N[j']}\\
&&\lesssim \sum_{k'<k-10}\left\|P_{j'}\left[\partial^{\alpha}U^N_{k'}\partial_{\alpha}\left(u_{<k-10}\ud_k-u_k\ud_{<k-10}\right)\phi\right]\right\|_{N[j']}\\
&&\lesssim\sum_{k'<k-10}2^{-\kappa(j-k')}c_{k'}c_k\|\phi\|_{S[j]}\\
&&\lesssim 2^{-\kappa(j-k)}\epsilon \,c_k\|\phi\|_{S[j]}.
\end{eqnarray*}
Let us now consider the term $$P_{j'}\left[\left(U^N_{<k-10}\partial^{\alpha}u_{<k-10}\partial_{\alpha}u_k^{\dagger}\right)\phi\right],$$
from term $IV$.
We have, by (\ref{eq:boundforUk}) and the trilinear estimate, 
\begin{eqnarray*}
&&\left\|P_{j'}\left[\left(U^N_{<k-10}\partial^{\alpha}u_{<k-10}\partial_{\alpha}u_k^{\dagger}\right)\phi\right]\right\|_{N[j']}\\
&&\lesssim \sum_{k'<k-10}2^{-\kappa(j-k')}c_k c_{k'}\|\phi\|_{S[j]}\lesssim 2^{-\kappa (j-k)}\epsilon \,c_k\|\phi\|_{S[j]},
\end{eqnarray*}
for $j'=j+O(1)$.\\
The term $$P_{j'}\left[\left(U^N_{<k-10}\partial^{\alpha}u_k\partial_{\alpha}\ud_{<k-10}\right)\phi\right]$$ can be controlled similarly.\\
It remains to control term $III$. For this, we need to use the equation for $u$. Since $u$ satisfies the wave map equation, we see that
\begin{equation}\label{eq:equforuk'}
\Box u_{k'}=P_{k'}\left(u\,\partial^{\alpha}\ud\partial_{\alpha}u\right),\,\,{\rm for\,\,each\,\,}k'\leq k.
\end{equation}
It suffices to show that, for any $\varphi$ with Fourier support $2^{j-3}\leq |\xi|\leq 2^{j+3}$ and $k'<j-6$,
\begin{equation}\label{eq:estimatethird}
\left\|P_{j'}\left[P_{k'}\left(u\,\partial^{\alpha}\ud\partial_{\alpha}u\right)\varphi\right]\right\|_{N[j']}\lesssim 2^{-\nu(j-k')}\epsilon^{\frac{1}{2}} c_{k'}\|\varphi\|_{S[j]},
\end{equation}
for $j'=j+O(1)$. Indeed, from (\ref{eq:estimatethird}), it follows that
\begin{eqnarray*}
&&\left\|P_{j'}\left[U^N_{<k-10}\Box u_{<k-10}\ud_k\phi\right]\right\|_{N[j']}\\
&&\lesssim \sum_{k'<k-10}\left\|P_{j'}\left[U^N_{<k-10}\Box u_{k'}\ud_k\phi\right]\right\|_{N[j']}\\
&&\lesssim \sum_{k'<k-10}2^{-\nu (j-k')}c_{k'}c_k\|\phi\|_{S[j]}\\
&&\lesssim \epsilon\,c_k 2^{-\nu (j-k)}\|\phi\|_{S[j]}
\end{eqnarray*}
and that
\begin{eqnarray*}
&&\left\|P_{j'}\left[U^N_{<k-10}\Box u_{k}\,\ud_{<k-10}\phi\right]\right\|_{N[j']}\\
&&\lesssim 2^{-\nu (j-k)}c_k\,\epsilon^{\frac{1}{2}}\|\phi\|_{S[j]}.
\end{eqnarray*}
These estimates are sufficient for the completion of the induction, due to the presence of the extra $\epsilon^{\frac{1}{2}}$ factor, which can be used to absorb various constants in the inequalities.\\
To prove (\ref{eq:estimatethird}), let us decompose $P_{k'}\left(u\,\partial^{\alpha}\ud\partial_{\alpha}u\right)\phi$ as
\begin{eqnarray*}
&&P_{k'}\left(u\,\partial^{\alpha}\ud\partial_{\alpha}u\right)\varphi=P_{k'}\left(u_{>\frac{j+k'}{2}}\,\partial^{\alpha}\ud\partial_{\alpha}u\right)\varphi\\
&&\hspace{.4in}+\,P_{k'}\left(u_{\leq \frac{j+k'}{2}}\,\partial^{\alpha}\ud\partial_{\alpha}u\right)\varphi=I_1+I_2.
\end{eqnarray*}
For $I_1$, by the trilinear estimates and symmetry, we can estimate as follows
\begin{eqnarray*}
&&\left\|P_{j'}\left[P_{k'}\left(u_{>\frac{j+k'}{2}}\,\partial^{\alpha}\ud\partial_{\alpha}u\right)\varphi\right]\right\|_{N[j']}=\left\|\sum_{k_2,\,k_3,\,k_1>\frac{j+k'}{2}}P_{j'}\left[P_{k'}\left(u_{k_1}\,\partial^{\alpha}\ud_{k_2}\partial_{\alpha}u_{k_3}\right)\varphi\right]\right\|_{N[j']}\\
&&\lesssim \|\varphi\|_{S[j]}\left(\sum_{k_1>\frac{j+k'}{2},\,k_3\ge k_1+O(1),\,k_3=k_2+O(1)}2^{-\kappa(\max_{1\leq i\leq 3}k_i-k')}2^{-\kappa (k_1-\min\{k_2,\,k_3\})_+}c_{k_1}c_{k_2}c_{k_3}\right.\\
&&\quad\quad\quad+\left.\sum_{k_1>\frac{j+k'}{2},\,k_2<k_3-C,\,k_3= k_1+O(1)}2^{-\kappa(\max_{1\leq i\leq 3}k_i-k')}2^{-\kappa (k_1-\min\{k_2,\,k_3\})_+}c_{k_1}c_{k_2}c_{k_3}\right)\\
&&\lesssim\epsilon^2\, \|\varphi\|_{S[j]}\left(\sum_{k_1>\frac{j+k'}{2}}c_{k_1}2^{-\kappa(k_1-k')}+\sum_{k_1>\frac{j+k'}{2}}c_{k_1}2^{-\kappa (k_1-k')}\right)\\
&&\lesssim \epsilon^2\,2^{-\frac{\kappa}{2}(j-k')} \|\varphi\|_{S[j]}c_{k'}.
\end{eqnarray*}
Now let us deal with the term $I_2=P_{k'}\left(u_{\leq \frac{j+k'}{2}}\,\partial^{\alpha}\ud\partial_{\alpha}u\right)\varphi$. In this case, we can insert $P_{<\frac{j+k'}{2}+C}$ in front of $\partial^{\alpha}\ud\,\partial_{\alpha}u$, use symmetry, and obtain that
\begin{eqnarray*}
&&\left\|P_{j'}\left[P_{k'}\left(u_{\leq \frac{j+k'}{2}}\,\partial^{\alpha}\ud\partial_{\alpha}u\right)\varphi\right]\right\|_{N[j']}\\
&&=\left\|P_{j'}\left[P_{k'}\left(u_{\leq \frac{j+k'}{2}}\,P_{<\frac{j+k'}{2}+C}\left(\partial^{\alpha}\ud\partial_{\alpha}u\right)\right)\varphi\right]\right\|_{N[j']}\\
&&\lesssim\left\|\sum_{k_1\leq k_2,\,k_2=k_1+O(1)}P_{j'}\left[P_{k'}\left(u_{\leq \frac{j+k'}{2}}\,P_{<\frac{j+k'}{2}+C}\left(\partial^{\alpha}\ud_{k_1}\partial_{\alpha}u_{k_2}\right)\right)\varphi\right]\right\|_{N[j']}+\\
&&\hspace{.4in}+\,\left\|\sum_{k_1\leq k_2-C}P_{j'}\left[P_{k'}\left(u_{\leq \frac{j+k'}{2}}\,P_{<\frac{j+k'}{2}+C}\left(\partial^{\alpha}\ud_{k_1}\partial_{\alpha}u_{k_2}\right)\right)\varphi\right]\right\|_{N[j']},
\end{eqnarray*}
which can be estimated as
\begin{eqnarray*}
&&\lesssim \left\|\sum_{k_1\leq k_2,\,k_2=k_1+O(1),\,k_1>\frac{3j+k'}{4}}P_{j'}\left[P_{k'}\left(u_{\leq \frac{j+k'}{2}}\,P_{<\frac{j+k'}{2}+C}\left(\partial^{\alpha}\ud_{k_1}\partial_{\alpha}u_{k_2}\right)\right)\varphi\right]\right\|_{N[j']}\\
&&\hspace{.3in}+\,\left\|\sum_{k_1\leq k_2,\,k_2=k_1+O(1),\,k_1\leq \frac{3j+k'}{4}}P_{j'}\left[P_{k'}\left(u_{\leq \frac{j+k'}{2}}\,P_{<\frac{j+k'}{2}+C}\left(\partial^{\alpha}\ud_{k_1}\partial_{\alpha}u_{k_2}\right)\right)\varphi\right]\right\|_{N[j']}\\
&&\hspace{.3in}+\,\sum_{k_1\leq k_2-C,\,k_2\leq \frac{j+k'}{2}+C}2^{-\kappa (j-k_1)}\|\varphi\|_{S[j]}c_{k_1}c_{k_2}
\end{eqnarray*}
which is
\begin{eqnarray*}
&&\lesssim \sum_{k_1>\frac{3j+k'}{4}}2^{-\kappa\left(k_1-\frac{j+k'}{2}\right)}\|\varphi\|_{S[j]}\cdot c_{k_1}^2+\sum_{k_1\leq\frac{3j+k'}{4}}c_{k_1}^2\cdot 2^{-\kappa(j-k_1)}\|\varphi\|_{S[j]}\\
&&\hspace{.4in}+\,\sum_{k_1\leq k_2-C,\,k_2\leq \frac{j+k'}{2}+C}2^{-\kappa(j-k_1)}\|\varphi\|_{S[j]}c_{k_1}c_{k_2}\\
&&\lesssim 2^{-\frac{\kappa}{8}(j-k')}\epsilon\,c_{k'}\|\varphi\|_{S[j]}.
\end{eqnarray*}
Combining the above estimates for $I,\,II,\,III,\,IV$ terms, the claim follows.\\

\medskip
\noindent
{\it Substep (2): Control of the term containing $\partial \psi$.}\\
Now we address the main term in the nonlinearity that forced us to use the gauge transform
$$\wt{h}=\left[U_{<m-10}\left(u_{<m-10}\partial^{\alpha}\ud_{<m-10}-\partial^{\alpha}u_{<m-10}\ud_{<m-10}\right)-\partial^{\alpha}U_{<m-10}\right]\partial_{\alpha}\psi.$$
Note that by (\ref{eq:iterativeU}), we have
\begin{eqnarray*}
&&-\wt{h}=\left[\partial^{\alpha}U_{<m-10}-U_{<m-10}\left(u_{<m-10}\partial^{\alpha}\ud_{<m-10}-\partial^{\alpha}u_{<m-10}\ud_{<m-10}\right)\right]\partial_{\alpha}\psi\\
&&=\sum_{k<m-10}\left[\partial^{\alpha}U_k-U_{<m-10}\left(u_{<m-10}\partial^{\alpha}\ud_k-\partial^{\alpha}u_k\ud_{<m-10}\right)\right]\partial_{\alpha}\psi\\
&&=\sum_{k<m-10}\left[\partial^{\alpha}U_k-U_{<k-10}\left(u_{<k-10}\partial^{\alpha}\ud_k-\partial^{\alpha}u_k\ud_{<k-10}\right)\right]\partial_{\alpha}\psi\\
&&\hspace{.5in}-\,\sum_{k<m-10}U_{k-10\leq \cdot<m-10}\left(u_{<m-10}\partial^{\alpha}\ud_k-\partial^{\alpha}u_k\,\ud_{<m-10}\right)\partial_{\alpha}\psi\\
&&\hspace{.5in}-\,\sum_{k<m-10}U_{<k-10}\left(u_{k-10\leq \cdot<m-10}\partial^{\alpha}\ud_k-\partial^{\alpha}u_k\,\ud_{k-10\leq\cdot<m-10}\right)\partial_{\alpha}\psi\\
&&=\sum_{k<m-10}\left[\pua U_{<k-10}\left(u_{<k-10}\ud_k-u_k\ud_{<k-10}\right)\right.\\
&&\hspace{.8in}\left.-\,U_{<k-10}\left(\partial^{\alpha}u_{<k-10}\ud_k-u_k\partial^{\alpha}\ud_{<k-10}\right)\right]\partial_{\alpha}\psi+\mathcal{R}.
\end{eqnarray*}
To estimate the $\mathcal{R}$ term, let us firstly bound for $m'=m+O(1)$,
\begin{eqnarray*}
&&\left\|P_{m'}\left[U_{k-10\leq \cdot<m-10}u_{<m-10}\partial^{\alpha}\ud_k\partial_{\alpha}\psi\right]\right\|_{N[m']}\\
&&\lesssim \sum_{k-10\leq k'<m-10}\left\|P_{m'}\left[U_{k'}u_{<m-10}\partial^{\alpha}\ud_k\partial_{\alpha}\psi\right]\right\|_{N[m']}\\
&&\lesssim  \sum_{k-10\leq k'<m-10}\sup_{m''=m+O(1)}\left\|P_{m''}\left[U_{k'}\partial^{\alpha}\ud_k\partial_{\alpha}\psi\right]\right\|_{N[m'']}\\
&&\lesssim \sum_{k-10\leq k'<m-10}2^{-\kappa (k'-k)}c_{k'}c_kc_m\lesssim \sum_{k-10\leq k'<m-10}2^{-(\kappa-\vartheta)(k'-k)}c_k^2c_m\\
&&\lesssim c_k^2c_m.
\end{eqnarray*}
Other terms in $\mathcal{R}$ can be treated similarly. Thus
\begin{equation}
\sup_{m'=m+O(1)}\left\|P_{m'}\left[\mathcal{R}\right]\right\|_{N[m']}\lesssim \sum_kc_k^2c_m\lesssim \epsilon^2c_m,
\end{equation}
and consequently $\mathcal{R}$ is disposable.\\
We can estimate for $m'=m+O(1)$
\begin{eqnarray*}
&&\left\|P_{m'}\left[\partial^{\alpha}U_{<k-10}\left(u_{<k-10}\ud_k-u_k\ud_{<k-10}\right)\partial_{\alpha}\psi\right]\right\|_{N[m']}\\
&&\lesssim \sum_{k'<k-10}\left\|P_{m'}\left[\partial^{\alpha}U_{k'}\left(u_{<k-10}\ud_k-u_k\ud_{<k-10}\right)\partial_{\alpha}\psi\right]\right\|_{N[m']}\\
&&\lesssim \sum_{k'<k-10}2^{-\kappa (k-k')}c_kc_{k'}\|\psi\|_{S[m]}\\
&&\lesssim \sum_{k'<k-10}2^{-(\kappa-\vartheta)(k-k')}c_k^2\|\psi\|_{S[m]}\lesssim c_k^2\|\psi\|_{S[m]},
\end{eqnarray*}
and
\begin{eqnarray*}
&&\left\|P_{m'}\left[U_{<k-10}\partial^{\alpha}u_{<k-10}\ud_k\partial_{\alpha}\psi\right]\right\|_{N[m']}\\
&&\lesssim  \sum_{k'<k-10}\left\|P_{m'}\left[U_{<k-10}\partial^{\alpha}u_{k'}\ud_k\partial_{\alpha}\psi\right]\right\|_{N[m']}\\
&&\lesssim \sum_{k'<k-10}2^{-\kappa (k-k')}c_{k'}c_k\|\psi\|_{S[m]}\\
&&\lesssim c_k^2\|\psi\|_{S[m]}.
\end{eqnarray*}
Thus in summary, we can estimate
\begin{equation*}
\sup_{m'=m+O(1)}\|\wt{h}\|_{N[m']}\lesssim \sum_{k<m-10}c_k^2\|\psi\|_{S[m]}\lesssim \epsilon^2c_m,
\end{equation*}
and consequently $\wt{h}$ is disposable.\\

\noindent
{\it Substep (3): $U_{<m-10}f$ term is disposable.}\\
This follows directly as $f$ is disposable.\\

\medskip
\noindent
{\bf Step 5: Proof of the channel of energy inequality for the good frequency piece.}\\
Take $m\in\mathcal{K}$. By the estimates from {\it Step 4}, we can write the equation for $w$ in Step 3 as 
\begin{equation}\label{eq:finalequforw}
\partial_{tt}w-\Delta w=h,
\end{equation}
with $h$ being disposable in the generalized sense, that is, $h=\lim\limits_{k\to\infty}h_k$ in the sense of distributions and $\sup\limits_{m'=m+O(1)}\|P_{m'}h_k\|_{N[m']}\lesssim \epsilon c_m$ uniformly in k. Let us now study how the outgoing condition (\ref{eq:outgoingfork}) on the initial data of $\psi$ has been transformed. Recall that 
$$w=U_{<m-10}\psi.$$
Hence $$\nabla_{x,t}w=\nabla_{x,t}U_{<m-10}\psi+U_{<m-10}\nabla_{x,t}\psi.$$
Thus at time $t=0$, by (\ref{eq:boundforUk}) and the outgoing condition for $\psi$,
\begin{eqnarray*}
&&\|\nabla_{x,t}w(0)\|_{L_x^2(B^c_{1+\lambda}\cup B_{1-\lambda})}\\
&&\lesssim \|\nabla_{x,t}U_{<m-10}(0)\|_{L^{2}}\|\psi(0)\|_{L_x^{\infty}}+\|U_{<m-10}(0)\|_{L_x^{\infty}}\|\nabla_{x,t}\psi(0)\|_{L^2(B_{1+\lambda}^c\cup B_{1-\lambda})}\\
&&\lesssim \left(\sum_{k<m-10}\|\nabla_{x,t}U_k(0)\|^2_{L^2_x}\right)^{\frac{1}{2}}\|P_m(u_0,\,u_1)\|_{\HL}+\delta^{\frac{1}{100}}\|P_m(u_0,\,u_1)\|_{\HL}\\
&&\lesssim \left(\sum_{k<m-10}c_k^2\right)^{\frac{1}{2}}\|P_m(u_0,\,u_1)\|_{\HL}+\delta^{\frac{1}{100}}\|P_m(u_0,\,u_1)\|_{\HL}\\
&&\lesssim \left(\epsilon+\delta^{\frac{1}{100}}\right)\|P_m(u_0,\,u_1)\|_{\HL}.
\end{eqnarray*}
Similar calculations show that 
$$\|\spartial w_0\|_{L^2}+\|\partial_rw_0+w_1\|_{L^2}\lesssim\left(\epsilon+\delta^{\frac{1}{100}}\right)\|P_m(u_0,\,u_1)\|_{\HL},$$
and
$$\|(w_0,\,w_1)\|_{\HL}\ge (1-\gamma(\epsilon))\|P_m(u_0,\,u_1)\|_{\HL},$$
with a suitable $\gamma\to 0$ as $\epsilon\to 0$.
If $\delta$ and $\epsilon$ are chosen sufficiently small, then by the channel of energy inequality for the linear wave equation and the bound on $h$, we conclude using (\ref{eq:inhomogeneousenergyestimates}) and (\ref{eq:Senergy}) that for all $t\ge 0$, 
\begin{equation}\label{eq:channelforw}
\int_{|x|\ge \frac{\beta+1}{2}+t}|\nabla_{x,t}w|^2(x,t)\,dx\ge \left|\frac{3+\beta}{4}\right|\|P_m(u_0,\,u_1)\|^2_{\HL}-C\epsilon^2c_m^2.
\end{equation}
Since $\nabla_{x,t}U_{<m-10}\,\psi$ is small in $L^2$ (smaller than $C\epsilon c_m$) and $U_{<m-10}$ is almost orthogonal by (\ref{eq:orthSc}), the channel of energy inequality (\ref{eq:channelform}) for $\psi$ follows (again by choosing $\epsilon,\,\delta$ sufficiently small depending on $\beta$). \\
This finishes the proof of the theorem.\\

It remains to prove Claim \ref{claim:Ijareperturbative} and Claim \ref{claim:commutingperturbative}.\\

\noindent
{\it Proof of Claim \ref{claim:Ijareperturbative}:} We need to control $I_1,\,I_2,\,I_3$.\\
For $I_1$, by the trilinear estimate and symmetry, we get that for $m'=m+O(1)$
\begin{eqnarray*}
&&\left\|P_{m'}\left(u_{\ge m-10}\partial^{\alpha}\ud\partial_{\alpha}u\right)\right\|_{N[m']}\\
&&=\left\|\sum_{k_1\ge m-10,\,k_2,\,k_3}P_{m'}\left(u_{k_1}\partial^{\alpha}\ud_{k_2}\partial_{\alpha}u_{k_3}\right)\right\|_{N[m']}\\
&&\lesssim \sum_{k_1\ge m-10,\,k_2\ge k_3}2^{-\kappa\left(\max\{k_1,\,k_2,\,k_3\}-m\right)_+}2^{-\kappa(k_1-\min\{k_2,\,k_3\})_+}\times\\
&&\hspace{1.2in}\times\,\|u_{k_1}\|_{S[k_1]}\,\|u_{k_2}\|_{S[k_2]}\,\|u_{k_3}\|_{S[k_3]}\\
&&\\
&&\lesssim \sum_{k_1\ge m-10,\,k_2\ge k_3}2^{-\kappa\left(\max\{k_1,\,k_2\}-m\right)_+}2^{-\kappa (k_1-k_3)_+}c_{k_1}\epsilon^2\\
&&\lesssim \epsilon^2\,c_m\sum_{k_1\ge m-10,\,k_2\ge k_3}2^{-(\kappa-\vartheta)\left(\max\{k_1,\,k_2\}-m\right)_+}2^{-\kappa (k_1-k_3)_+}\lesssim \epsilon^2\,c_m.
\end{eqnarray*}
 For $I_2$, by the product property and null form estimate, we get that for $m'=m+O(1)$
\begin{eqnarray*}
&&\left\|P_{m'}\left(u_{<m-10}\partial^{\alpha}\ud_{>m+10}\partial_{\alpha}u\right)\right\|_{N[m']}\\
&&\lesssim \sum_{k_1>m+10,\,k_2=k_1+O(1)}\left\|P_{m'}\left(u_{<m-10}\partial^{\alpha}\ud_{k_1}\partial_{\alpha}u_{k_2}\right)\right\|_{N[m']}\\
&&\lesssim \sum_{k_1>m+10,\,k_2=k_1+O(1)}2^{-\kappa (k_1-m)}\|u_{k_1}\|_{S[k_1]}\|u_{k_2}\|_{S[k_2]}\\
&&\lesssim  \sum_{k_1>m+10,\,k_2=k_1+O(1)}2^{-\kappa (k_1-m)}c_{k_1}c_{k_2}\\
&&\lesssim \epsilon\,c_m\sum_{k_1>m+10,\,k_2=k_1+O(1)}2^{-(\kappa-\vartheta) (k_1-m)}\lesssim \epsilon\,c_m.
\end{eqnarray*}
For $I_3$, by the product property and null form estimate, we get that for $m'=m+O(1)$
\begin{eqnarray*}
&&\left\|P_{m'}\left(u_{<m-10}\partial^{\alpha}\ud_{m-10\leq\cdot\leq m+10}\partial_{\alpha}u_{\ge m-10}\right)\right\|_{N[m']}\\
&&\lesssim \sum_{k\ge m-10}\left\|P_{m'}\left(u_{<m-10}\partial^{\alpha}\ud_{m-10\leq\cdot\leq m+10}\partial_{\alpha}u_{k}\right)\right\|_{N[m']}\\
&&\lesssim \sum_{k\ge m-10}2^{-\kappa (k-m)}\epsilon\,c_m\lesssim \epsilon\,c_m.
\end{eqnarray*}
Thus the terms $I_1,\,I_2,\,I_3$ 
are all disposable. The claim is proved.\\

\medskip
\noindent
{\it Proof of Claim \ref{claim:commutingperturbative}:} Noting that $$P_m\left(u_{m-10\leq\cdot\leq m+10}\right)=P_mu=\psi,$$ by Lemma \ref{lm:commutinglemma}, we get that 
\begin{eqnarray*}
&&P_m\left(u_{<m-10}\partial_{\alpha}\ud_{<m-10}\partial^{\alpha}u_{m-10\leq\cdot\leq m+10}\right)-u_{<m-10}\partial_{\alpha}\ud_{<m-10}\partial^{\alpha}\psi\\
&&=2^{-m}L\left(\nabla \left(u_{<m-10}\partial_{\alpha}\ud_{<m-10}\right),\,\partial^{\alpha}u_{m-10\leq\cdot\leq m+10}\right)\\
&&=2^{-m}L\left(\nabla u_{<m-10}\partial_{\alpha}\ud_{<m-10},\,\partial^{\alpha}u_{m-10<\cdot<m+10}\right)+\\
&&\hspace{.3in}+\,2^{-m}L\left(u_{<m-10}\partial_{\alpha}\nabla\ud_{<m-10},\,\partial^{\alpha}u_{m-10<\cdot<m+10}\right).
\end{eqnarray*}
Thus, noting that
$$\|\nabla u_k\|_{S[k]}\lesssim 2^k\|u_k\|_{S[k]},$$
 by the trilinear estimate for the first term in the above and the product estimate and null form estimate for the second, we get that for $m'=m+O(1)$
\begin{eqnarray*}
&&\left\|P_{m'}\left[P_m\left(u_{<m-10}\partial_{\alpha}\ud_{<m-10}\partial^{\alpha}u_{m-10<\cdot<m+10}\right)-u_{<m-10}\partial_{\alpha}\ud_{<m-10}\partial^{\alpha}\psi\right]\right\|_{N[m']}\\
&&\lesssim \sum_{k_1<m-10,\,k_2<m-10}2^{-m}\left\|P_{m'}\left[L\left(\nabla u_{k_1}\partial_{\alpha}\ud_{k_2},\,\partial^{\alpha}u_{m-10<\cdot<m+10}\right)\right]\right\|_{N[m']}\\
&&\hspace{.3in}+\,\sum_{k<m-10}2^{-m}\left\|P_{m'}\left[L\left(u_{<m-10}\partial_{\alpha}\nabla\ud_{k},\,\partial^{\alpha}u_{m-10<\cdot<m+10}\right)\right]\right\|_{N[m']}\\
&&\lesssim 2^{-m}\sum_{k_1<m-10,\,k_2<m-10}2^{k_1}2^{-\kappa(k_1-k_2)_+}c_m\epsilon^2+\sum_{k<m-10}2^{-m}2^k\epsilon \,c_m\lesssim \epsilon c_m.
\end{eqnarray*}
The first part of the claim is proved. The proof of the second part is similar.

\end{section}

\begin{section}{Morawetz estimates and applications}
In the previous sections, the main tools we use are all perturbative in nature. In order to understand the dynamics of large wave maps, we need some global control on the solution. Such global control is often achieved with help of suitable monotonicity formulae. The most important monotonicity formula here are the energy flux identity and the Morawetz estimate. This section follows similar arguments in Sterbenz-Tataru \cite{Tataru4}. We use less geometric, but perhaps slightly more transparent notations similar to those used in \cite{DJKM}. \\

For notational convenience, we shall work with a classical wave map $u$ defined on $R^2\times(0,1]$, that equals $u_{\infty}\in S^2$ for large $x$. Let us firstly look at the energy flux identity for $u$. We thus have
$$\partial_{tt}u-\Delta u=(|\nabla u|^2-|u_t|^2)u,\,\,{\rm in}\,\,R^2\times(0,1].$$
Noting that $\ud\cdot u_t\equiv 0$, we have the identity
\begin{equation}\label{eq:utidentity}
\left(\partial_{tt}\ud-\Delta \ud\right)\cdot u_t=0,\,\,{\rm in}\,\,R^2\times(0,1].
\end{equation}
 Take $0<t_1<t_2<1$, and integrate the identity (\ref{eq:utidentity}) in the truncated lightcone $\{(x,t):\,|x|<t,\,t_1<t<t_2\}$, we obtain that
\begin{eqnarray*}
&&\int_{|x|<t_2}\nablaxtu(x,t_2)\,dx-\int_{|x|<t_1}\nablaxtu(x,t_1)\,dx\\
&&\hspace{.3in}-\frac{1}{\sqrt{2}}\int_{t_1}^{t_2}\int_{|x|=t}\flux\,d\sigma dt=0.
\end{eqnarray*}
Denote 
\begin{equation*}
\Fl(t_1,t_2):=\frac{1}{\sqrt{2}}\int_{t_1}^{t_2}\int_{|x|=t}\flux\,d\sigma dt
\end{equation*}
as the ``energy flux" through the lateral boundary of the lightcone. We see that
\begin{eqnarray*}
\Fl(t_1,t_2)&=&\int_{|x|<t_2}\nablaxtu(x,t_2)\,dx\\
&&\hspace{.4in}-\int_{|x|<t_1}\nablaxtu(x,t_1)\,dx.
\end{eqnarray*}
Since $\Fl(t_1,t_2)\ge0$, it follows that 
\begin{equation*}
\int_{|x|<t}\nablaxtu (x,t)\,dx
\end{equation*}
is nondecreasing, and has a limit as $t\to 0+$. Thus $\Fl(t_1,t_2)\to 0+$ as $t_1,\,t_2\to 0+$. \\

The control of energy flux plays an essential role in the following Morawetz estimate.
\begin{theorem}\label{th:Morawetz}
Let $u$ be a classical wave map with energy $E$ on $R^2\times (0,1]$ and $\epsilon\in(0,1)$. For each $0<\bt<1$, if $\Fl(0,\bt)<\epsilon E$, then
\begin{equation}\label{eq:Morawetz}
\int_{\epsilon \bt}^{\bt}\int_{|x|<t}\rho_{\epsilon\bt}^3\left(X^{\alpha}\partial_{\alpha}u\right)^2\,dxdt+\int_{|x|<\bt}\bt\rho_{\epsilon\bt}\left(\frac{|\nabla u|^2}{2}+\frac{|u_t|^2}{2}+\frac{x}{\bt}\cdot\nabla \ud\,u_t\right)(x,\bt)\,dx\lesssim E,
\end{equation}
where we set $\rho_{\epsilon\bt}:=\left((t+\epsilon\bt)^2-|x|^2\right)^{-\frac{1}{2}}$ and $X^{\alpha}=x^{\alpha}$ if $\alpha=1,2$, $X^{0}=t+\epsilon\bt$.
\end{theorem}

\smallskip
\noindent
{\it Proof.} By rescaling, we can assume without loss of generality that $\bt=1$. (Then $u$ is rescaled to $R^2\times \left(0,\,\frac{1}{\bt}\right]$) Let us integrate the identity 
$$\partial^{\alpha}\partial_{\alpha}\ud\,\rho_{\epsilon}\,X^{\beta}\partial_{\beta}u=0$$
on $\{(x,t):\,|x|<t,\,\epsilon<t<1\}$. We have
\begin{eqnarray*}
0&=&\int_{\epsilon}^1\int_{|x|<t}\partial^{\alpha}\partial_{\alpha}\ud\,\rho_{\epsilon}X^{\beta}\partial_{\beta}u\,dxdt\\
&=&\int_{\epsilon}^1\int_{|x|<t}\rho_{\epsilon}\,X^{\beta}\partial^{\alpha}(\partial_{\alpha}\ud\,\partial_{\beta}u)-\rho_{\epsilon}\,X^{\beta}\partial_{\beta}\frac{\partial^{\alpha}\ud\partial_{\alpha}u}{2}\,dxdt\\
&=&B+I,
\end{eqnarray*}
where the boundary term $B$ and the interior term $I$ are
\begin{eqnarray*}
B&=&\int_{\epsilon}^1\int_{|x|=t}\rho_{\epsilon}\,X^{\beta}\,n^{\alpha}\,\partial_{\alpha}\ud\,\partial_{\beta}u-\rho_{\epsilon}\,X^{\beta}\,n_{\beta}\frac{\partial^{\alpha}\ud\partial_{\alpha}u}{2}\,d\sigma dt\\
&&\,-\,\int_{|x|<1}\rho_{\epsilon}\,X^{\beta}\,\partial_t\ud\,\partial_{\beta}u(x,1)\,dx-\int_{|x|<1}(1+\epsilon)\rho_{\epsilon}\,\frac{\partial^{\alpha}\ud\partial_{\alpha}u}{2}(x,1)\,dx\\
&&\,+\,\int_{|x|<\epsilon}\rho_{\epsilon}\,X^{\beta}\,\partial_t\ud\,\partial_{\beta}u(x,\epsilon)\,dx+\int_{|x|<\epsilon}2\epsilon\,\rho_{\epsilon}\,\frac{\partial^{\alpha}\ud\partial_{\alpha}u}{2}(x,\epsilon)\,dx;
\end{eqnarray*}
and
\begin{equation*}
I=-\int_{\epsilon}^1\int_{|x|<t}\partial^{\alpha}\left(\rho_{\epsilon}\,X^{\beta}\right)\,\partial_{\alpha}\ud\,\partial_{\beta}u-\partial_{\beta}\left(\rho_{\epsilon}\,X^{\beta}\right)\,\frac{\partial^{\alpha}\ud\partial_{\alpha}u}{2}\,dxdt.
\end{equation*}
In the above we use the notation $n=\frac{1}{\sqrt{2}}\left(\frac{x}{|x|},\,-1\right)$, $n^j=n_j=\frac{1}{\sqrt{2}}\frac{x_j}{|x|}$ for $j=1,2$ and $n^0=-n_0=\frac{1}{\sqrt{2}}$. Hence $X^{\beta}n_{\beta}=-\frac{\epsilon}{\sqrt{2}}$ on $|x|=t$. We can compute 
$$\partial_j\rho_{\epsilon}=x_j\rho_{\epsilon}^3,\,\,\,\partial_t\rho_{\epsilon}=-t\rho_{\epsilon}^3-\epsilon\,\rho_{\epsilon}^3.$$ 
Hence $$X^{\beta}\partial_{\beta}\rho_{\epsilon}=-\rho_{\epsilon}.$$
We also note that $\epsilon\,\rho\leq 1$ when $|x|<t$, and record the following simple bound when $|x|=t$
$$\left|\rho_{\epsilon}\right|= \left(2\epsilon \,t+\epsilon^2\right)^{-\frac{1}{2}}\leq \epsilon^{-\frac{1}{2}}t^{-\frac{1}{2}}.$$
We can simplify the $B,\,I$ terms as
\begin{eqnarray*}
B&=&\frac{1}{\sqrt{2}}\int_{\epsilon}^1\int_{|x|=t}\rho_{\epsilon}\,\left(X^{\beta}\partial_{\beta}u\right)\cdot\frac{\left(x^{\alpha}\partial_{\alpha}u\right)}{t}+\epsilon\,\rho_{\epsilon}\frac{\partial^{\alpha}u^{\dagger}\,\partial_{\alpha}u}{2}\,d\sigma dt\\
&&\,\,\,\,-\int_{|x|<1}\rho_{\epsilon}\,\fluxOne(x,1)\,dx+O(E),
\end{eqnarray*}
and
\begin{eqnarray*}
-I&=&\int_{\epsilon}^1\int_{|x|<t}\bigg[X^{\beta}\partial_{\beta}\ud\,\partial^{\alpha}\rho_{\epsilon}\,\partial_{\alpha}u+\rho_{\epsilon}\,\left(|\nabla u|^2-|\partial_tu|^2\right)\\
&&\hspace{.4in}\left.-\,\frac{3}{2}\rho_{\epsilon}\,\left(\partial^{\alpha}\ud\,\partial_{\alpha}u\right)-X^{\beta}\partial_{\beta}\rho_{\epsilon}\,\frac{\partial^{\alpha}\ud\,\partial_{\alpha}u}{2}\right]\,dxdt\\
&=&\int_{\epsilon}^1\int_{|x|<t}\rho_{\epsilon}^3\left|X^{\beta}\partial_{\beta}u\right|^2\,dxdt.
\end{eqnarray*}
We can estimate
\begin{eqnarray*}
&&\left|\frac{1}{\sqrt{2}}\int_{\epsilon}^1\int_{|x|=t}\rho_{\epsilon}\,\left(X^{\beta}\partial_{\beta}u\right)\cdot\frac{\left(x^{\alpha}\partial_{\alpha}u\right)}{t}+\epsilon\,\rho_{\epsilon}\frac{\partial^{\alpha}u^{\dagger}\,\partial_{\alpha}u}{2}\,d\sigma dt\right|\\
&&\lesssim \int_{\epsilon}^1\int_{|x|=t}\rho_{\epsilon}\,t\left|\partial_tu+\frac{x}{t}\cdot\nabla u\right|^2\,d\sigma dt+\\
&&\hspace{.3in}+\,\int_{\epsilon}^1\int_{|x|=t}\epsilon\,\rho_{\epsilon}\,\left(\frac{|\nabla u|^2}{2}+\frac{|\partial_tu|^2}{2}+\frac{x}{t}\cdot\nabla \ud\,\partial_tu\right)\,d\sigma dt\\
&&\lesssim \epsilon^{-\frac{1}{2}}\Fl(0,1)\leq \epsilon^{\frac{1}{2}}E.
\end{eqnarray*}
Hence, combining the $B$ and $I$ terms, we conclude that 
\begin{eqnarray*}
&&\int_{|x|<1}\rho_{\epsilon}\,\fluxOne(x,1)\,dx+\int_{\epsilon}^1\int_{|x|<1}\rho_{\epsilon}^3\left|X^{\beta}\partial_{\beta}u\right|^2\,dxdt\\
&&\lesssim E.
\end{eqnarray*}
The theorem is proved.\\

Theorem \ref{th:Morawetz} has the following corollary.
\begin{corollary}\label{cor:vanishing}
Let $u$ be as above. For any $\tau_n\to 0+$, $\gamma_n\to 1-$ as $n\to\infty$, we have that
\begin{equation}\label{eq:corvani}
\int_{B_{\tau_n}\backslash B_{\gamma_n\tau_n}}\flux(x,\tau_n)\,dx=o_n(1).
\end{equation}
\end{corollary}

\smallskip
\noindent
{\it Proof:} Let $\epsilon_n:=2\,\Fl(0,\tau_n)/E$, then $\epsilon_n\to 0$ as $n\to\infty$. Theorem \ref{th:Morawetz} implies that
\begin{equation*}
\int_{B_{\tau_n}}\tau_n\,\rho_{\epsilon_n\tau_n}\flux(x,\tau_n)\,dx\lesssim E.
\end{equation*}
Note that 
$$\tau_n\,\rho_{\epsilon_n\tau_n}\gtrsim \left((1+\epsilon_n)^2-\gamma_n^2\right)^{-\frac{1}{2}}\to\infty,$$
for $t=\tau_n,\,|x|\in (\gamma_n\tau_n,\,\tau_n)$, we conclude that (\ref{eq:corvani}) holds.

\end{section}

\begin{section}{Universal blow up profile along a sequence of times}
Our goal in this section is to prove Theorem \ref{th:smallsolitonresolutionMain} along a sequence of times. Again for the ease of notations, we shall consider classical wave map $u$ defined on $R^2\times (0,1]$ that blows up at time $t=0$. Recall from (\ref{eq:concentrationradiusast}) the definition of $r(\epsilon_{\ast},\,t)$. By the small data theory and finite speed of propagation, we have $\lim\limits_{t\to0+} r(\epsilon_{\ast},\, t)=0$. That is, energy concentrates in smaller and smaller regions as $t\to 0+$. By the definition of $r(\epsilon_{\ast},\,t)$, we can find $x_{\ast}(t)$ such that
\begin{equation}
\|\OR{u}(t)\|_{\HL\left(B_{2r(\epsilon_{\ast},\,t)}(x_{\ast}(t))\right)}>\frac{\epsilon_{\ast}}{2},
\end{equation}
for $t$ close to $0$. Again by small data global existence and finite speed of propagation,  $x_{\ast}(t)$ remains in a bounded region for $t\in (0,\,1]$. Assume without loss of generality that $x_{\ast}(t_n)\to 0$ as $t_n\to 0$ along a sequence of times $t_n$. Since $r(\epsilon_{\ast},\,t)\to 0$, we see that for any $r>0$,
\begin{equation}\label{eq:singularconcentration}
\liminf_{t\to 0}\left\|\OR{u}(t)\right\|_{\HL(B_r)}>\frac{\epsilon_{\ast}}{2}.
\end{equation}
In general, we call a point $\overline{x}$ {\it singular} if for any $r>0$
$$\limsup_{t\to 0}\|\OR{u}(t)\|_{\HL(B_r(\overline{x}))}>\frac{\epsilon_{\ast}}{2}.$$
By finite speed of propagation and energy flux identity, this is equivalent to requiring that for any $r>0$,
$$\liminf_{t\to 0}\|\OR{u}(t)\|_{\HL(B_r(\overline{x}))}>\frac{\epsilon_{\ast}}{2}.$$
 As the energy is conserved and finite, there can only be finitely many singular points, and in particular the singular points are isolated. \\

For the singular point $x_{\ast}=0$, since singular points are isolated, there exists $r_1>0$ such that for any $\overline{x}\in B_{r_1}\backslash\{0\}$, $\bx$ is not a singular point. Hence we can find $\tr>0$ with 
\begin{equation}\label{eq:small1}
\|\OR{u}(\tau_n)\|_{\HL( B_{\tr}(\overline{x}))}<\epsilon_{\ast}
\end{equation}
along a sequence of times $\tau_n\to 0$. In particular, we have
$$\|\spartial u(\cdot,\,\tau_n)\|_{L^2(B_{\tr}(\bx))}<\epsilon_{\ast}.$$
Hence there exists $\br_n\in \left(\frac{\tr}{2},\,\tr\right)$ with
$$\int_{\partial B_{\br_n}(\bx)}|\spartial u|^2(x,\tau_n)\,d\sigma \lesssim \frac{\epsilon_{\ast}^2}{\br_n}.$$
Denoting $\overline{u}_n$ as the average of $u(\tau_n)$ over $\partial B_{\br_n}(\bx)$, that is
$$\overline{u}_n=\frac{1}{2\pi \br_n}\int_{\partial B_{\br_n}(\bx)}u(\tau_n)\,d\sigma.$$
By Sobolev inequality, we get that
\begin{equation}\label{eq:small0}
\left\|u(\tau_n)-\overline{u}_n\right\|_{L^{\infty}(|x-\bx|=\br_n)}\lesssim \epsilon_{\ast}.
\end{equation}
Take smooth cutoff function $\eta_n\in C_c^{\infty}\left(B_{2\br_n}(\bx)\right)$ with $\eta_n|_{B_{\br_n}(\bx)}\equiv 1$ and $|\nabla \eta_n|\lesssim (\br_n)^{-1}$. Recall that for any $v\in R^2$ with $v\neq 0$, 
$$Pv=\frac{v}{|v|}.$$
Define
\begin{equation*}
(u_{0n},\,u_{1n})=\left\{\begin{array}{lr}
             (u,\,\partial_tu)(\tau_n) & {\rm in}\,\,B_{\br_n}(\bx);\\
      \left(P[\eta_n (u(\br_n \theta,\tau_n)-\overline{u}_n)+\overline{u}_n],\,0\right) & {\rm in}\,\,\left(B_{\br_n}(\bx)\right)^c.
                                     \end{array}\right.
\end{equation*}
 By (\ref{eq:small1}) and (\ref{eq:small0}), direct computation shows that
\begin{equation}\label{eq:small3}
\|(u_{0n},\,u_{1n})\|_{\HL}\lesssim \epsilon_{\ast}.
\end{equation}
Hence by small energy global existence theory and finite speed of propagation, we see that the solution $u_n$ to the wave map equation with $\OR{u}_n(\tau_n)=(u_{0n},\,u_{1n})$ is global and that 
\begin{equation}\label{eq:localequal}
u_n\equiv u\,\,\,{\rm for}\,\,\, |x-\bx|<\frac{\br_n}{4}\,\,{\rm and}\,\,\, t\in(0,\,\tau_n]
\end{equation}
 for sufficiently large $n$. Since $u_n\in C\left([0,\,\tau_n],\,\HL\right)$ and (\ref{eq:localequal}) holds, we conclude that $u$ can be extended to $t=0$ so that $u\in C\left([0,\,\tau_n],\,\HL(B_{\tr/8}(\bx))\right)$. Since $\bx\in B_{r_{1}}\backslash\{0\}$ is arbitrary, we conclude that $u$ can be extended to $t=0$ in $B_{r_1}$ with $u\in C\left([0,\,1],\,\HL(B_{r_1}\backslash B_r)\right)$ for each $0<r<r_1$.\\
In addition, by the regularity of $u_n$, we also have the additional (qualitative) regularity condition that $u\in C^{\infty}\left(B_{r_1}\times [0,\,1]\backslash\{(0,0)\}\right)$. One can of course apply the same argument to other singular points. As a result, we see that $u\in C^{\infty}(R^2\times[0,1]\backslash \{(x_j,\,0)\})$ where $x_j$ are the singular points.\\
On the other hand, since $\OR{u}(t)$ is bounded in $\HL$ and $|u|\equiv 1$, we can extract a weak limit $(v_0,\,v_1)\in\HL$ along a sequence of times $t_n\to 0+$. This limit is in fact a strong limit outside an arbitrarily small neighborhood of the finitely many singular points. From the above analysis, $(v_0,\,v_1)\in C^{\infty}(R^2\backslash\{x_j\})$. Let 
$$v=u$$
for $\inf\limits_j|x-x_j|>t$. Then $v\in C^{\infty}\left(R^2\times [0,1]\backslash \bigcup\limits_j \{|x-x_j|\leq t,\,t\in[0,\,1]\}\right)$, and by the same arguments as in the proof of Lemma \ref{lm:exteriornice},
\begin{equation}\label{eq:regularpart}
\lim_{t\to 0+}\bigg\|\OR{v}(\cdot,\,t)-(v_0,\,v_1)\bigg\|_{\HL\left(\bigcap\limits_j\{|x-x_j|>t\}\right)}=0.
\end{equation}
We shall call $v$ the {\it regular part} of the wave map $u$. The main issue is to understand the behavior of the wave map $u$ inside singularity lightcones $ \bigcup\limits_j \{|x-x_j|\leq t,\,t\in[0,\,1]\}$.\\

We shall prove
\begin{theorem}\label{th:solitonsequential}
Let $u$ be a classical wave map with energy 
\begin{equation}\label{eq:energyconstraint5}
\E(\OR{u})<\E(Q,0)+\epsilon_0^2,
\end{equation}
 where $Q$ is a harmonic map with degree $1$, defined on $R^2\times (0,\,1]$ that blows up at time $t=0$ with the origin being a singular point. Assume that $\epsilon_0$ is sufficiently small. Then there exists a sequence of times $t_n\to 0+$,  $\ell\in R^2$ with $|\ell|\ll 1$, $x_n\in R^2,\,\lambda_n>0$ with 
$$\lim_{n\to\infty}\frac{x_n}{t_n}=\ell,\,\,\,\lambda_n=o(t_n),$$
and $(v_0,\,v_1)\in \HL\cap C^{\infty}(R^2\backslash\{0\})$, such that
\begin{equation}
\OR{u}(t_n)=(v_0,\,v_1)+\left(Q_{\ell},\,\lambda_n^{-1}\partial_tQ_{\ell}\right)\left(\frac{x-x_n}{\lambda_n},\,\frac{t-t_n}{\lambda_n}\right)+o_{\HL}(1),
\end{equation}
as $n\to\infty$.
\end{theorem}

\smallskip
\noindent
{\it Remark:} As we discussed in the introduction, the main new point in Theorem \ref{th:solitonsequential} is that we eliminate any possible energy concentration near the boundary of singularity lightcone $|x|<t$. By the energy constraint (\ref{eq:energyconstraint5}) and the proof below, there is only one singularity point. Hence $u$ is regular outside $\{(x,t):\,|x|\leq t\}$. The main task is to understand the behavior of $u$ inside $\{(x,t):\,|x|\leq t\}$.

\smallskip
\noindent
{\it Proof:} Our starting point is the work of Grinis \cite{Grinis}, which completely characterized the concentration of energy in $\{(x,\,\tau_n):\,|x|<a\tau_n\}$ for any $a\in(0,1)$ as traveling waves, for a suitable time sequence $\tau_n\to 0+$. See Theorem 1.1 and Theorem 1.2 in \cite{Grinis}.  In our case, due to the energy constraint (\ref{eq:energyconstraint5}), there can only be one traveling wave. Hence, as a particular consequence of a rescaled version of the {\it asymptotic decomposition} in Theorem 1.2 of \cite{Grinis}, we have for $|x|<\tau_n$,
\begin{equation}\label{eq:Grinis}
\OR{u}(\tau_n):=\left(Q_{\ell}\left(\frac{x-x_n}{r_n},\,0\right),\,r_n^{-1}\nabla \partial_tQ_{\ell}\left(\frac{x-x_n}{r_n},\,0\right)\right)+(w_{0n},\,w_{1n})+o_{\HL}(1),
\end{equation}
as $n\to\infty$, where $|\ell|\ll 1$, $r_n=o(\tau_n)$, $\ell=\lim\limits_{n\to\infty}\frac{x_n}{\tau_n}$ and
\begin{equation}\label{eq:residuenull}
\int_{|x|<a\tau_n}|\nabla w_{0n}|^2+|w_{1n}|^2\,dx\to 0,
\end{equation}
as $n\to\infty$ for any $a\in(0,1)$. Our main task is to show that 
$$\int_{|x|<\tau_n}|\nabla w_{0n}|^2+|w_{1n}|^2\,dx\to 0$$
as $n\to\infty$. By (\ref{eq:residuenull}), we have to prove that that for any $\gamma_n\to 1-$,
\begin{equation}\label{eq:noboundaryenergy3}
\limsup_{n\to\infty}\int_{B_{\tau_n}\backslash B_{\gamma_n\tau_n}}\nablaxtu(x,\tau_n)\,dx=0,
\end{equation}
assuming that 
\begin{equation}
\label{eq:wn0}
\limsup_{n\to\infty}\int_{B_{\gamma_n\tau_n}}|\nabla w_{0n}|^2+|w_{1n}|^2\,dx=0. 
\end{equation} 

We now apply the channel of energy inequality and prove (\ref{eq:noboundaryenergy3}). 

Suppose that (\ref{eq:noboundaryenergy3}) is not true. Then there exists $\epsilon_2>0$, such that, by passing to a subsequence if necessary, we have for all sufficiently large $n$,
\begin{equation}\label{eq:concentrationepsilon2}
\E_n^2:=\int_{B_{\tau_n}\backslash B_{\gamma_n\tau_n}}\nablaxtu(x,\tau_n)\,dx\ge \epsilon^2_2.
\end{equation}
By the energy constraint, we must also have
\begin{equation}\label{eq:energyonboundaryissmall}
\int_{B_{\frac{\tau_n}{2}}^c\cap B_{\tau_n}}\nablaxtu(x,\tau_n)\,dx\lesssim\epsilon_0^2.
\end{equation}
Corollary \ref{cor:vanishing} implies that
\begin{equation}\label{eq:morawetzoutgoing}
\int_{B_{\tau_n}\backslash B_{\gamma_n\tau_n}}\left(\frac{|\nabla u|^2}{2}+\frac{|\partial_tu|^2}{2}+\partial_t\ud\,\partial_r u\right)(x,\tau_n)\,dx=o_n(1).
\end{equation}
Since $u$ is regular for $|x|>t$, we have for any $r>0$,
\begin{equation}\label{eq:smallexterior17}
\limsup_{t\to 0+}\int_{B_{2r}\backslash B_t}\nablaxtu(x,t)\,dx\leq \delta(r)\to 0,\,\,{\rm as}\,\,r\to 0+.
\end{equation}
Fix a small $r>0$ whose value is to be determined below. We can find $r_{1n}\in \left(\frac{r}{2},\,r\right)$, $r_{2n}\in \left(\frac{\tau_n}{2},\,\frac{3}{4}\tau_n\right)$, such that
$$\int_{\partial B_{r_{1n}}}|\spartial u|^2(\tau_n)\,d\sigma\lesssim \frac{\delta(r)}{r_{1n}},\,\,\,{\rm and}\,\,\,\int_{\partial B_{r_{2n}}}|\spartial u|^2(\tau_n)\,d\sigma=\frac{o_n(1)}{r_{2n}}.$$
Let
$$\overline{u}^1_n=\frac{1}{2\pi r_{1n}}\int_{\partial B_{r_{1n}}}u(\tau_n)\,d\sigma,\,\,\,{\rm and}\,\,\,\overline{u}^2_n=\frac{1}{2\pi r_{2n}}\int_{\partial B_{r_{2n}}}u(\tau_n)\,d\sigma.$$
Fix radial $\eta_{1n}\in C_c^{\infty}(B_{2r_{1n}})$ with $\eta_{1n}|_{B_{r_{1n}}}\equiv 1$, and radial $1-\eta_{2n}\in C_c^{\infty}(B_{r_{2n}})$ with $1-\eta_{2n}|_{B_{\frac{r_{2n}}{2}}}\equiv 1$. Define
\begin{equation}
(u_{0n},\,u_{1n})=\left\{\begin{array}{lr}
     \left(P\left[\eta_{1n}\,\left(u(r_{1n}\theta,\,\tau_n)-\overline{u}^1_n\right)+\overline{u}^1_n\right],\,0\right)&{\rm in}\,\,B_{r_{1n}}^c;\\
                \OR{u}(\tau_n)&{\rm in}\,\,B_{r_{1n}}\backslash B_{r_{2n}};\\
     \left(P\left[\eta_{2n}\,\left(u(r_{2n}\theta,\,\tau_n)-\overline{u}^2_n\right)+\overline{u}^2_n\right],\,0\right)&{\rm in}\,\,B_{r_{2n}}
                                    \end{array}\right. 
\end{equation}
Then for sufficiently large $n$, in view of (\ref{eq:wn0}) and (\ref{eq:smallexterior17}),
$$\|(u_{0n},\,u_{1n})\|_{\HL\left(B_{\tau_n}^c\bigcup B_{\tau_n\gamma_n}\right)}\lesssim \delta(r),$$
and
$$\epsilon_0\gtrsim \|(u_{0n},\,u_{1n})\|_{\HL}>\E_n+O(\delta(r)) \geq \epsilon_2+O(\delta(r)).$$
In addition, by (\ref{eq:morawetzoutgoing}), for sufficiently large $n$
$$\|u_{1n}+\partial_ru_{0n}\|_{L^2}+\|\spartial u_{0n}\|_{L^2}\lesssim \delta(r).$$
Let $u_n$ be the solution to the wave map equation with $\OR{u}(\tau_n)=(u_{0n},\,u_{1n})$. Then if $r$ is taken sufficiently small so that $\delta(r)$ is much smaller than $\epsilon_2$ by (a rescaled and time translated version of) Theorem \ref{th:channelofenergywavemap} we conclude that for $t\ge \tau_n$
\begin{equation}\label{eq:thinchannel}
\int_{|x|>t-\frac{\tau_n}{8}}\left|\nabla_{x,t}u_n\right|^2(x,t)\,dx\gtrsim \E_n^2.
\end{equation}
Take $t=\frac{r}{8}$ in (\ref{eq:thinchannel}), we get that for all sufficiently large $n$,
\begin{equation}\label{eq:thinchannel'}
\int_{|x|>\frac{r}{8}-\frac{\tau_n}{8}}\left|\nabla_{x,t}u_n\right|^2(x,\frac{r}{8})\,dx\gtrsim \E_n^2.
\end{equation}
By the energy inequality, (\ref{eq:smallexterior17}) and the definition of $u_n$, we see that for $t\leq \frac{r}{8}$,
\begin{equation}\label{eq:noleakage2}
\int_{|x|>t}\left|\nabla_{x,t}u_n\right|^2(x,t)\,dx\lesssim\delta(r).
\end{equation}
By finite speed of propagation, we also have $u\equiv u_n$ for $t-\frac{\tau_n}{4}<|x|<\frac{r}{4}$ and $t\leq \frac{r}{8}$. Combining with (\ref{eq:thinchannel'}), we conclude that
\begin{equation}\label{eq:thinchannel''}
\int_{\frac{r}{8}>|x|>\frac{r}{8}-\frac{\tau_n}{4}}\left|\nabla_{x,t}u\right|^2(x,\frac{r}{8})\,dx\gtrsim \E^2_n\ge\epsilon_2^2>0,
\end{equation}
if we choose $r$ sufficiently small, so that $\delta(r)$ is much smaller than $\epsilon_2^2$. However, (\ref{eq:thinchannel''}) contradicts with the fact that $\OR{u}\left(\frac{r}{8}\right)\in\HL$ for sufficiently large $n$. \\

Therefore, combining the above with the regular part outside the singularity lightcone, we get that along the sequence $\tau_n$, 
\begin{equation}\label{eq:resolution1}
\OR{u}(\tau_n)=(v_0,\,v_1)+(Q_{\ell},\,r_n^{-1}\partial_tQ_{\ell})\left(\frac{x-x_n}{r_n},\,0\right)+o_{\HL}(1),\,\,\,{\rm as}\,\,n\to\infty.
\end{equation}
The theorem is proved.
\end{section}

\begin{section}{Coercivity and universal profile for all times}

Our next task is to use a rigidity property of the energy to extend the decomposition we obtained from the last section to all times. One important tool is the following coercivity property of the energy functional near traveling waves. 
\begin{theorem}\label{th:coercivity}
Let $\mathcal{M}_1$ be the space of harmonic maps from $R^2$ to $S^2$ with topological degree $1$. Fix $\ell\in R^2$ with $|\ell|<1$ and for any $Q\in\mathcal{M}_1$, let $Q_{\ell}$ be the Lorentz transform of $Q$ with velocity $\ell$, that is
\begin{equation}\label{eq:Ql}
Q_{\ell}(x,t)=Q\left(x-\frac{\ell\cdot x}{|\ell|^2}\ell+\frac{\frac{\ell\cdot x}{|\ell|^2}\ell-\ell t}{\sqrt{1-|\ell|^2}}\right).
\end{equation}
Denote $\mathcal{M}_{1,\ell}$ as the space of $Q_{\ell}$ with $Q\in\mathcal{M}_1$. For $0<\epsilon<\epsilon_0$ and $\epsilon_0$ sufficiently small,
suppose that $(v_0,\,v_1)\in\HL$, with $|v_0(x)|\equiv 1$ and $v_0^{\dagger}\cdot v_1\equiv 0$, satisfies
\begin{eqnarray}
&&\deg(v_0)=1;\label{eq:degree1}\\
&&\left|\int_{R^2}\partial_{x_j}v_0^{\dagger}\,v_1\,dx-\int_{R^2}\partial_{x_j}Q_{\ell}^{\dagger}\,\partial_tQ_{\ell}\,dx\right|<\epsilon;\label{eq:equalmomentum}\\
&&\int_{R^2}\left(\frac{|\nabla v_0|^2}{2}+\frac{|v_1|^2}{2}\right)dx\leq \int_{R^2}\left(\frac{|\partial_tQ_{\ell}|^2}{2}+\frac{|\nabla Q_{\ell}|^2}{2}\right)\,dx+\epsilon;\label{eq:equalenergy}\\
&&\inf_{Q\in\mathcal{M}_1}\left\|(v_0,\,v_1)-(Q_{\ell},\,\partial_tQ_{\ell})\right\|_{\HL}<\epsilon_0.\label{eq:energyconstraint23}
\end{eqnarray}
Then there exists $\delta(\epsilon)>0$ with $\delta(\epsilon)\to 0$ as $\epsilon\to 0$, such that
\begin{equation}\label{eq:closetotravelingwave}
\inf_{Q\in\mathcal{M}_1}\left\|(v_0,\,v_1)-(Q_{\ell},\,\partial_tQ_{\ell})\right\|_{\HL}<\delta(\epsilon).
\end{equation}

\end{theorem}
\smallskip
\noindent
{\it Remark:} As we will see in the proof 
\begin{eqnarray*}
&&\int_{R^2}\frac{|\nabla Q_{\ell}|^2}{2}+\frac{|\partial_tQ_{\ell}|^2}{2}\,dx=\frac{4\pi}{\sqrt{1-|\ell|^2}},\\
&&-\int_{R^2}\partial_tQ^{\dagger}_{\ell}\,\partial_{x_j}Q_{\ell}\,dx=\frac{4\pi\ell_j}{\sqrt{1-|\ell|^2}},
\end{eqnarray*}
for any $Q\in\mathcal{M}_1$. The conditions (\ref{eq:equalmomentum}) and (\ref{eq:equalenergy}) are thus independent of the choice of $Q\in \mathcal{M}_1$.

The definition of degree $\deg(f)$ for mappings between manifolds is classical. For the definition with $f:\,R^2\to S^2\subset R^3$ and $f\in \dot{H}^1$ used here, we refer to \cite{Brezis}, see in particular (1) in page 205 of \cite{Brezis}. We also remark that the harmonic maps in $\mathcal{M}_1$ have been completely characterized as degree $1$ rational functions (M\"obius transforms), see \cite{Elles} and a more recent discussion in \cite{Oh}.  By elementary geometric properties of M\"obius transforms,  it is easy to see that degree one harmonic maps from $R^2$ to $S^2\subset R^3$ are unique up to the symmetries of $R^2$ and $S^2$. More precisely in an appropriate coordinate system, the harmonic maps in $\mathcal{M}_1$ are co-rotational. 

\smallskip
\noindent
{\it Proof.} Without loss of generality, let us assume that $\ell=le_1=(l,\,0)$.  Suppose that (\ref{eq:closetotravelingwave}) is false, then for each $n=1,2,\dots$, by symmetry, we can assume that there exist $(v_{0n},\,v_{1n})\in\HL$ with $|v_{0n}|\equiv 1$ and $v_{0n}^{\dagger}\,v_{1n}\equiv 0$, such that 
\begin{eqnarray}
&&\deg(v_{0n})=1;\label{eq:degree1'}\\
&&\left|\int_{R^2}\partial_{x_j}v_{0n}^{\dagger}v_{1n}\,dx-\int_{R^2}\partial_{x_j}Q_{\ell}^{\dagger}\,\partial_tQ_{\ell}\,dx\right|<\frac{1}{n};\label{eq:equalmomentum'}\\
&&\int_{R^2}\left(\frac{|\nabla v_{0n}|^2}{2}+\frac{|v_{1n}|^2}{2}\right)dx\leq \int_{R^2}\left(\frac{|\partial_tQ_{\ell}|^2}{2}+\frac{|\nabla Q_{\ell}|^2}{2}\right)\,dx+\frac{1}{n};\label{eq:equalenergy'}\\
&&\inf_{Q\in\mathcal{M}_1}\left\|(v_{0n},\,v_{1n})-(Q_{\ell},\,\partial_tQ_{\ell})\right\|_{\HL}<\epsilon_0.\label{eq:energyconstraint23'}
\end{eqnarray}
In addition, 
\begin{equation}\label{eq:closetotravelingwave'}
\inf_{Q\in\mathcal{M}_1}\left\|(v_{0n},\,v_{1n})-(Q_{\ell},\,\partial_tQ_{\ell})\right\|_{\HL}>\delta_0>0.
\end{equation}
For fixed $(v_0,\,v_1)\in\HL$, with $|v_0(x)|\equiv 1$ and $v_0^{\dagger}\cdot v_1\equiv 0$, assume without loss of generality that $v_0$ is positively oriented, that is, 
\begin{equation}\label{eq:degreeH1}
\deg(v_0)=-\frac{1}{4\pi}\int_{R^2}v_0^{\dagger}\cdot\left(\partial_{1}v_0\times\partial_2v_0\right).
\end{equation}

Consider the following algebraic identity
\begin{multline}\label{eq:magicidentity}
\int_{R^2}\left(\frac{|\nabla v_0|^2}{2}+\frac{|v_1|^2}{2}\right)\,dx\\
=\frac{1}{2}\int_{R^2}\left|v_{1}+l\partial_1v_{0}\right|^2\,dx+\frac{1}{4}\int_{R^2}\left|\sqrt{1-l^2}\partial_1v_{0}-v_{0}\times\partial_2v_0\right|^2\,dx\\
+\frac{1}{4}\int_{R^2}\left|\partial_2v_0+\sqrt{1-l^2}v_0\times\partial_1v_0\right|^2\,dx-\sqrt{1-l^2}\int_{R^2}v_0^{\dagger}\cdot\left(\partial_1v_0\times\partial_2v_0\right)\,dx\\
-l\int_{R^2}\partial_1v_0^{\dagger}\,v_1\,dx.
\end{multline}
(\ref{eq:magicidentity}) is a modified form of the remarkable decomposition of energy in \cite{Belavin}, see also the illuminating discussion in page 3 of \cite{RodSter}. The modification here is necessary in order to take into account the momentum part.\\
Direct calculations show that
\begin{equation}\label{eq:infoQl}
\int_{R^2}\left(\frac{|\partial_tQ_{\ell}|^2}{2}+\frac{|\nabla Q_{\ell}|^2}{2}\right)\,dx=\frac{4\pi}{\sqrt{1-l^2}}\,\,\,{\rm and}\,\,\,-\int_{R^2}\partial_tQ_{\ell}^{\dagger}\,\partial_1Q_{\ell}\,dx=\frac{4\pi l}{\sqrt{1-l^2}}.
\end{equation}
We can assume, after rotation, that $(v_{0n},\,v_{1n})$ has the same momentum as $(Q_{\ell_n},\,\partial_tQ_{\ell_n})$ with $\ell_n=l_n\,e_1$. Then $|l_n-l|\lesssim \frac{1}{n}$.
Applying (\ref{eq:magicidentity}) to $(v_{0n},\,v_{1n})$ and using the assumptions on $(v_{0n},\,v_{1n})$, we get that
\begin{eqnarray*}
&&\int_{R^2}\left(\frac{|\nabla v_{0n}|^2}{2}+\frac{|v_{1n}|^2}{2}\right)\,dx\\
&&=\frac{1}{2}\int_{R^2}\left|v_{1n}+l_n\partial_1v_{0n}\right|^2\,dx+\frac{1}{4}\int_{R^2}\left|\sqrt{1-l_n^2}\partial_1v_{0n}-v_{0n}\times\partial_2v_{0n}\right|^2\,dx\\
&&\hspace{.3in}+\,\frac{1}{4}\int_{R^2}\left|\partial_2v_{0n}+\sqrt{1-l_n^2}v_{0n}\times\partial_1v_{0n}\right|^2\,dx-\sqrt{1-l_n^2}\int_{R^2}v_{0n}^{\dagger}\cdot\left(\partial_1v_{0n}\times\partial_2v_{0n}\right)\,dx\\
&&\hspace{.5in}-\,l_n\int_{R^2}\partial_1v_{0n}^{\dagger}\,v_{1n}\,dx\\
&&=\frac{1}{2}\int_{R^2}\left|v_{1n}+l\partial_1v_{0n}\right|^2\,dx+\frac{1}{4}\int_{R^2}\left|\sqrt{1-l^2}\partial_1v_{0n}-v_{0n}\times\partial_2v_{0n}\right|^2\,dx\\
&&\hspace{.3in}+\,\frac{1}{4}\int_{R^2}\left|\partial_2v_{0n}+\sqrt{1-l^2}v_{0n}\times\partial_1v_{0n}\right|^2\,dx+\frac{4\pi}{\sqrt{1-l^2}}+O\left(\frac{1}{n}\right).
\end{eqnarray*}
In the above we used the expression for degree and momentum.
From (\ref{eq:equalenergy'}) and (\ref{eq:infoQl}), we conclude that
\begin{multline}\label{eq:magicsmall}
\frac{1}{2}\int_{R^2}\left|v_{1n}+l\partial_1v_{0n}\right|^2\,dx+\frac{1}{4}\int_{R^2}\left|\sqrt{1-l^2}\partial_1v_{0n}-v_{0n}\times\partial_2v_{0n}\right|^2\,dx\\
\hspace{.3in}+\,\frac{1}{4}\int_{R^2}\left|\partial_2v_{0n}+\sqrt{1-l^2}v_{0n}\times\partial_1v_{0n}\right|^2\,dx\\
=O\left(\frac{1}{n}\right).
\end{multline}
By (\ref{eq:energyconstraint23'}), applying suitable symmetry transformation to $(v_{0n},\,v_{1n})$ if necessary, we can assume that for suitable $\widetilde{Q}_{\ell}\in\mathcal{M}_{\ell,1}$, 
\begin{equation}
(v_{0n},\,v_{1n})=(\widetilde{Q}_{\ell},\,\partial_t\widetilde{Q}_{\ell})(x,0)+(r_{0n},\,r_{1n}),
\end{equation}
with $$\|(r_{0n},\,r_{1n})\|_{\HL}\leq 2\epsilon_0.$$
Passing to a subsequence, we can assume that $(v_{0n},\,v_{1n})\rightharpoonup (v_0,\,v_1)$ as $n\to\infty$,
with $$\left\|(v_0,\,v_1)-(\widetilde{Q}_{\ell},\,\partial_t\widetilde{Q}_{\ell})(x,0)\right\|_{\HL}\leq2\epsilon_0.$$
Hence by the continuity of topological degree \footnote{This is a direct consequence of the definition (\ref{eq:degreeH1}) of degree, and can be proved by noting that $\int_{R^2}\partial_xu \times \partial_yu\,dx\,dy=0$ for any $\dot{H}^1$ mapping from $R^2\to S^2$, and the dominated convergence theorem.} in $\dot{H}^1$  and the fact that degree only takes value in integers (see \cite{Brezis}), we see that if $\epsilon_0$ is taken small enough, then
$$\deg(v_0)=1.$$
(\ref{eq:magicsmall}) implies that $(v_0,\,v_1)$ satisfies the {\it first order} ``Bogomol'nyi equations" (see \cite{Bogo}):
\begin{eqnarray}
&&v_1+l\partial_1v_0=0;\nonumber\\
&&\sqrt{1-l^2}\partial_1v_0-v_0\times\partial_2v_0=0;\nonumber\\
&&\partial_2v_0+\sqrt{1-l^2}v_0\times\partial_1v_0=0.\label{eq:bogo}
\end{eqnarray}
Equations (\ref{eq:bogo}) can be reduced by an obvious change of variable to the case $l=0$, in which case they can be explicitly solved as harmonic maps. Hence we see that there exists $\widetilde{\widetilde{Q}}\in \mathcal{M}_1$ such that
$$(v_0,\,v_1)=\left(\widetilde{\widetilde{Q}}_{\ell},\,\partial_t\widetilde{\widetilde{Q}}_{\ell}\right)(x,0).$$
Thus we can write
$$(v_{0n},\,v_{1n})=\left(\widetilde{\widetilde{Q}}_{\ell},\,\partial_t\widetilde{\widetilde{Q}}_{\ell}\right)(x,0)+\left(\widetilde{r}_{0n},\,\widetilde{r}_{1n}\right),$$
with $\left(\widetilde{r}_{0n},\,\widetilde{r}_{1n}\right)\rightharpoonup 0$ as $n\to\infty$.
Then the energy expansion for $(v_{0n},\,v_{1n})$  around $\left(\widetilde{\widetilde{Q}}_{\ell},\,\partial_t\widetilde{\widetilde{Q}}_{\ell}\right)(x,0)$ gives
\begin{eqnarray*}
\E(Q_{\ell},\,\partial_t{Q}_{\ell})+\frac{1}{n}
&\ge&\int_{R^2}\left(\frac{|\nabla v_{0n}|^2}{2}+\frac{|v_{1n}|^2}{2}\right)\,dx\\
&=&{\displaystyle \frac{1}{2}\int_{R^2}\left|\nabla \widetilde{\widetilde{Q}}_{\ell}\right|^2+\left|\partial_t\widetilde{\widetilde{Q}}_{\ell}\right|^2\,dx}\\
&&\hspace{.3in}+\,\int_{R^2}\nabla\widetilde{\widetilde{Q}}_{\ell}^{\,\dagger}\,\nabla \widetilde{r}_{0n}+  \partial_t\widetilde{\widetilde{Q}}_{\ell}^{\,\dagger}\,\widetilde{r}_{1n}\,dx\\
&&\hspace{.3in}+\,\int_{R^2}\frac{|\nabla \widetilde{r}_{1n}|^2}{2}+\frac{|\widetilde{r}_{1n}|^2}{2}\,dx\\
&=&\E\left(\OR{\widetilde{\widetilde{Q}}_{\ell}}\right)+\,\int_{R^2}\frac{|\nabla \widetilde{r}_{1n}|^2}{2}+\frac{|\widetilde{r}_{1n}|^2}{2}\,dx+o_n(1).
\end{eqnarray*}
By (\ref{eq:equalenergy'}) and the fact that $\E\left(\OR{\widetilde{\widetilde{Q}}}_{\ell}\right)=\E(\OR{Q_{\ell}})$, we see that
$$\left(\widetilde{r}_{0n},\,\widetilde{r}_{1n}\right)\to 0,\,\,\,{\rm in}\,\,\HL.$$
This is a contradiction to (\ref{eq:closetotravelingwave'}). The Theorem is proved.\\

\medskip
Now we turn to the proof of the second main theorem in the paper. \\
\begin{theorem}\label{th:smallsolitonresolutionlast}
Let $u$ be a classical wave map defined on $R^2\times (0,\,1]$ with energy $\E(\OR{u})<\E(Q,0)+\epsilon_0^2$, where $Q$ is a harmonic map of degree $1$, that blows up at time $0$ and at the origin. Assume that $\epsilon_0$ is sufficiently small. Then there exist $\ell\in R^2$ with $|\ell|\ll 1$, $x(t)\in R^2,\,\lambda(t)>0$ with 
$$\lim_{t\to 0}\frac{x(t)}{t}=\ell,\,\,\,\lambda(t)=o\left(t\right),$$
and $(v_0,\,v_1)\in \HL\cap C^{\infty}(R^2\backslash\{0\})$ with $(v_0-u_{\infty},\,v_1)$ being compactly supported, such that
\begin{eqnarray*}
{\rm (i)}&&\inf\bigg\{\left\|\OR{u}(t)-(v_0,\,v_1)-\left(Q_{\ell},\,\partial_tQ_{\ell}\right)\right\|_{\HL}:\,Q_{\ell}\in \mathcal{M}_{\ell,1}\bigg\}\to 0,\,\,{\rm as}\,\,t\to 0;\\
&&\\
{\rm (ii)}&&\bigg\|\OR{u}(t)-(v_0,\,v_1)\bigg\|_{\HL\left(R^2\backslash B_{\lambda(t)}(x(t))\right)}\to 0\,\,{\rm as}\,\,t\to 0.
\end{eqnarray*}
\end{theorem}
\smallskip
\noindent
{\it Proof:} We have already proved that along a sequence of times $t_n\to 0+$, 
\begin{equation}\label{eq:predecomp}
\OR{u}(t_n)=\left(Q_{\ell}\left(\frac{x-x_n}{\lambda_n},0\right),\,\frac{1}{\lambda_n}\partial_tQ_{\ell}\left(\frac{x-x_n}{\lambda_n},0\right)\right)+(v_0,\,v_1)+o_{\HL}(1),
\end{equation}
where $(v_0,\,v_1)\in\HL\cap C^{\infty}(R^2\backslash\{0\})$ and
$$\lim_{n\to\infty}\frac{x_n}{t_n}=\ell,\,\,\,{\rm and}\,\,\,\lambda_n=o(t_n),\,\,\,{\rm as}\,\,n\to\infty.$$
Since
\begin{equation}\label{eq:sec6small}
\epsilon_n:=\int_{B_{2t_n}\backslash B_{t_n}}\nablaxtu(x,t_n)\,dx\to 0,\,\,\,{\rm as}\,\,n\to\infty,
\end{equation}
we can find $r_n\in(t_n,\,2t_n)$ such that
$$\int_{\partial B_{r_n}}\frac{|\spartial u|^2}{2}(x,t_n)\,d\sigma\lesssim\frac{\epsilon_n}{r_n}.$$
Let
$$\overline{u}_n=\frac{1}{2\pi r_n}\int_{\partial  B_{r_n}} u(t_n)\,d\sigma.$$

Take a smooth cutoff function $\eta_n$ with $\eta_n\equiv 1$ on $B_{r_n}$, ${\rm supp}\,\eta_n\Subset B_{2r_n}$ and $|\nabla\eta_n|\lesssim r_n^{-1}$.
Define
\begin{equation*}
(u_{0n},\,u_{1n})=\left\{\begin{array}{lr}
                       \OR{u}(x,t_n)&{\rm for}\,\,|x|<r_n;\\
                   \left(P\left[\eta_n(r)(u(r_n\theta,t_n)-\overline{u}_n)+\overline{u}_n\right],\,0\right)&{\rm for}\,\,|x|>r_n.
                                  \end{array}\right.
\end{equation*}
One can check that $(u_{0n},\,u_{1n})\in \dot{H}^s\times H^{s-1}$ for $s<\frac{3}{2}$, and $u_{0n}\equiv P(\overline{u}_n)$ for large $x$. Moreover,
\begin{equation}\label{eq:sec6small2}
\|(u_{0n},\,u_{1n})\|^2_{\HL(B^c_{t_n})}\lesssim \epsilon_n.
\end{equation}
Let $u_n$ be the solution to the wave map equation with $\OR{u}_n(t_n)=(u_{0n},\,u_{1n})$.\footnote{The local existence of $u_n$ follows from subcritical wellposedness theory. }  Then by finite speed of propagation, $u_n\equiv u$ for $|x|<t$ and $t\in(0,\,t_n]$, assuming that $u_n$ is defined in $[t,\,t_n]$. In addition, by (\ref{eq:sec6small2}) and energy flux identity, since the energy flux of $u_n$ is equal to that of $u$ on $|x|=t,\,t\in(0,\,t_n]$ which decays to zero as $n\to\infty$, we get that
\begin{equation}\label{eq:exteriorsmall}
\int_{|x|>t}|\nabla_{x,t}u_n|^2(x,t)\,dx\lesssim \epsilon_n+o_n(1),
\end{equation}
for $t\leq t_n$, again assuming that $u_n$ is defined in $[t,\,t_n]$. As $u_n$ is identical to $u$ in the singularity light cone $|x|<t,\,0<t\leq t_n$ and $u_n$ has small energy for $|x|\ge t,\,0<t\leq t_n$, we conclude that $u_n$ is defined for $t\in (0,\,t_n]$.
 From (\ref{eq:predecomp}), it is easy to verify that
\begin{eqnarray*}
&&{\rm deg}(u_n(t_n))=1;\\
&&\E(\OR{u}_n)\leq \E(\OR{Q}_{\ell})+o_n(1);\\
&&\left|\mathcal{M}(\OR{u}_n)-\mathcal{M}(\OR{Q}_{\ell})\right|=o_n(1),
\end{eqnarray*}
where $\mathcal{M}(\OR{u})$ denotes the momentum of $u$.
Hence by Theorem \ref{th:coercivity}, $\OR{u}_n(t)$ stays in a $\delta(\epsilon_n)$ neighborhood of $\mathcal{M}_{\ell,1}$ for $t\leq t_n$ with $\delta(\epsilon_n)\to 0$ as $n\to\infty$. It follows that
\begin{equation}\label{eq:convergencetosoliton}
\lim_{t\to 0}\,\inf\bigg\{\left\|\OR{u}(t)-(v_0,\,v_1)-\left(Q_{\ell},\,\partial_tQ_{\ell}\right)\right\|_{\HL}:\,Q_{\ell}\in \mathcal{M}_{\ell,1}\bigg\}=0.
\end{equation}
Part (i) of the theorem is proved. The fact that all degree $1$ harmonic maps are co-rotational implies that $\mathcal{M}_{\ell,1}$ is a compact set in the energy space, modulo translations and dilations. Hence, by the regularity of $u$ outside the singularity lightcone, we can find $x(t)$ and $\lambda(t)$ with $\lambda(t)=o(t)$ and 
$$\limsup_{t\to 0+}\frac{|x(t)|}{t}\leq 1.$$
The main remaining task is to show that
\begin{equation}\label{eq:thecenter}
\lim_{t\to 0+}\frac{x(t)}{t}=\ell.
\end{equation}
Without loss of generality, let us assume that $\ell=l\,e_1$
By (\ref{eq:convergencetosoliton}), it follows that
\begin{equation}\label{eq:momentuminside}
\lim_{t\to 0}\int_{|x|<t}\,-\partial_tu\,\partial_{x_1}u(x,t)\,dx=\frac{4l\pi}{\sqrt{1-l^2}}.
\end{equation}
Direct computation shows
\begin{eqnarray*}
&&\frac{d}{dt}\int_{|x|<t}x_1\nablaxtu(x,t)\,dx\\
&&\hspace{.3in}=\int_{|x|=t}x_1\nablaxtu(x,t)\,d\sigma+\int_{|x|=t}x_1\frac{x}{|x|}\cdot\nabla \ud \,\partial_tu\,d\sigma\\
&&\hspace{2in}-\int_{|x|<t}\partial_{x_1}u\,\partial_tu(x,t)\,dx.
\end{eqnarray*}
Integrating the above identity from $t=0$ to $t$, we get that
\begin{equation}\label{eq:centermotion}
\int_{|x|<t}x_1\nablaxtu(x,t)\,dx=O\left({\rm Flux}(0,t)\right)t+\frac{4l\pi\,t}{\sqrt{1-l^2}}+o(t).
\end{equation}
As ${\rm Flux}(0,t)\to 0$ as $t\to 0$, by (\ref{eq:centermotion}) and (\ref{eq:convergencetosoliton}), (\ref{eq:thecenter}) follows straightforwardly. The theorem is proved.
\end{section}

\end{document}